\documentclass{amsart}
\usepackage{amssymb,psfig,epsfig}

\title{On the Combinatorics of Crystal Graphs, I. Lusztig's Involution}

\author{Cristian Lenart}

\address{Department of Mathematics and Statistics, State University of New York, Albany, NY 12222}
\email{lenart@albany.edu}


\thanks{Cristian Lenart was supported by National Science Foundation 
grant DMS-0403029}


\numberwithin{equation}{section}

\theoremstyle{plain}
\newtheorem{theorem}{Theorem}[section]
\newtheorem{proposition}[theorem]{Proposition}
\newtheorem{lemma}[theorem]{Lemma}
\newtheorem{corollary}[theorem]{Corollary}

\theoremstyle{definition}
\newtheorem{definition}[theorem]{Definition}

\newtheorem{example}[theorem]{Example}

\theoremstyle{remark}
\newtheorem{remark}[theorem]{Remark}
\newtheorem{remarks}[theorem]{Remarks}

\def\R{\mathbb{R}}
\def\Z{\mathbb{Z}}
\def\C{\mathbb{C}}
\def\Q{\mathbb{Q}}

\def\O{\mathcal{O}}

\def\L{\mathcal{L}}

\def\Waff{W_{\mathrm{aff}}}

\def\hvee{h}
\def\h{\mathfrak{h}}
\def\hR{\mathfrak{h}^*_\mathbb{R}}
\newcommand{\casethree}[6]{\left\{ \begin{array}{ll} #1 &\mbox{if $#2$} \\[.05in]#3 &\mbox{if $#4$} \\[.05in] #5 & \mbox{if $#6$}\,. \end{array} \right.}

\newcommand{\casetwo}[3]{\left\{ \begin{array}{ll} #1 &\mbox{if $#2$} \\[.05in] #3 &\mbox{otherwise}\,. \end{array} \right.}

\newcommand{\casetwoex}[4]{\left\{ \begin{array}{ll} #1 &\mbox{if $#2$} \\[.05in] #3 &\mbox{if $#4$} \,. \end{array} \right.}
\newcommand{\casetwoexc}[4]{\left\{ \begin{array}{ll} #1 &\mbox{if $#2$} \\[.05in] #3 &\mbox{if $#4$} \,,\end{array} \right.}
\newcommand{\casethreeexc}[6]{\left\{ \begin{array}{ll} #1 &\mbox{if $#2$} \\[.05in] #3 &\mbox{if $#4$} \\[.05in] #5 &\mbox{if $#6$}\,,\end{array} \right.}

\newcommand{\floor}[1]{\lfloor #1 \rfloor}
\newcommand{\ceil}[1]{\lceil #1 \rceil}
\newcommand{\stacksum}[2]{\sum_{\begin{array}{c}\vspace{-5mm}\;\\ \vspace{-1mm}\scriptstyle{#1}\\ \scriptstyle{#2}\end{array}} }

\begin{document}
\bibliographystyle{plain}

\begin{abstract}
In this paper, we continue the development of a new combinatorial model
for the irreducible characters of a complex semisimple Lie group. This
model, which will be referred to as the alcove path model, can be viewed
as a discrete counterpart to the Littelmann path model. It leads to an
extensive generalization of the combinatorics of irreducible characters
from Lie type $A$ (where the combinatorics is based on Young tableaux, for
instance) to arbitrary type; our approach is type-independent. The main results of this paper are: (1) a
combinatorial description of the crystal graphs corresponding to the
irreducible representations (this result includes a transparent proof,
based on the Yang-Baxter equation, of the fact that the mentioned
description does not depend on the choice involved in our model); (2) a
combinatorial realization (which is the first direct generalization of {\em Sch\"utzenberger's involution} on tableaux) of Lusztig's involution on the
canonical basis exhibiting the crystals as self-dual
posets; (3) an analog for arbitrary root systems, based
on the Yang-Baxter equation, of Sch\"utzenberger's sliding algorithm, which is also known as jeu de taquin (this algorithm has many applications to the
representation theory of the Lie algebra of type $A$).
\end{abstract}

\maketitle

\section{Introduction}


We have recently given a simple combinatorial model for the irreducible characters of a complex semisimple Lie group $G$ and, more generally, 
for the {\em Demazure characters} \cite{LP}. For reasons explained below, we call our model the {\em alcove path model}. This was extended to complex symmetrizable Kac-Moody algebras in \cite{LP1} (that is, to infinite root systems). In this context, we also gave a {\em Littlewood-Richardson rule} for decomposing tensor products of irreducible representations and a branching rule. The exposition in \cite{LP} was in  the
context of the equivariant $K$-theory of the generalized flag variety $G/B$; more precisely, we first derived a {\em Chevalley-type multiplication formula} in $K_T(G/B)$, and then we deduced from it our Demazure character formula. By contrast, the exposition in \cite{LP1} was purely representation theoretic, being based on Stembridge's combinatorial model for Weyl characters \cite{St}. 

The alcove path model leads to an extensive generalization of the combinatorics of irreducible characters from Lie type $A$ (where the combinatorics is based on Young tableaux, for instance) to arbitrary type; our approach is type-independent. The present paper continues the study of the combinatorics of the new model, which was started in \cite{LP,LP1}. A future publication will be concerned with a direct generalization of the notion of the product of Young tableaux in the context of the product of crystals.

The main results of this paper are:
\begin{enumerate}
\item a combinatorial description of the crystal graphs corresponding to the
irreducible representations (Corollary \ref{all-cryst}); this result includes a transparent proof, based on the Yang-Baxter equation, of the fact that the mentioned
description does not depend on the choice involved in our model (Corollary \ref{isom-graphs});
\item a combinatorial realization of Lusztig's involution \cite{lusiqg} on the canonical basis (Theorem \ref{sch-inv-desc}, see also Example \ref{ex-inv}); this involution exhibits the crystals as self-dual posets, and corresponds to the action of the longest Weyl group element on an irreducible representation; our combinatorial realization is the first direct generalization of {\em Sch\"utzenberger's involution} on tableaux (see e.g. \cite{fulyt});
\item an analog for arbitrary root systems, based on the Yang-Baxter equation, of Sch\"utzenberger's {\em sliding} algorithm, which is also known as {\em jeu de taquin} (Section \ref{yb-moves}); this algorithm has many applications to the representation theory of the Lie algebra of type $A$ (see e.g. \cite{fulyt}). 
\end{enumerate}

Our model is based on the choice of an {\it alcove path}, which is a sequence of adjacent alcoves for
the {\em affine Weyl group} $\Waff$ of the Langland's dual group $G^\vee$. An alcove path  is best represented as a {\em $\lambda$-chain}, that is, as a sequence of positive roots corresponding to the common walls of successive alcoves in the mentioned sequence of alcoves. These chains extend the notion of a {\em reflection ordering} \cite{Dyer}. Given a fixed $\lambda$-chain, the objects that generalize semistandard Young tableaux are all the subsequences of roots that give rise to saturated increasing chains in Bruhat order (on the Weyl group $W$) upon multiplying on the right by the corresponding reflections. We call these subsequences {\em admissible subsets}. In \cite{LP1} we defined {\em root operators} on admissible subsets, which are certain partial operators associated with the simple roots; in type $A$, they correspond to the {\em coplactic operations} on tableaux \cite{lltpm}. The root operators produce a directed colored graph structure and a poset structure on admissible subsets. We showed in \cite{LP1} that this graph is isomorphic to the crystal graph of the corresponding irreducible representation if the chosen $\lambda$-chain is a special one. All this background information on the alcove path model is explained in more detail in Section \ref{alc-path}, following some general background material discussed in Section \ref{prelim}. 

In Section \ref{yb-moves}, we study certain discrete moves which allow us to deform any $\lambda$-chain into any other $\lambda$-chain (for a fixed dominant weight $\lambda$), and to biject the corresponding admissible subsets. We call these moves {\em Yang-Baxter moves} since they express the fact that certain operators satisfy the {\em Yang-Baxter equation}. We will explain below the reason for which the Yang-Baxter moves can be considered an analog of jeu de taquin for arbitrary root systems. We show that the Yang-Baxter moves commute with the root operators; this means that the directed colored graph defined by the root operators is invariant under Yang-Baxter moves, and it is thus independent from the choice of a $\lambda$-chain. Based on the special case in \cite{LP1} discussed above, this immediately implies that the mentioned graph is isomorphic to the corresponding crystal graph for any choice of a $\lambda$-chain. 

In Section \ref{sch-inv}, we present a combinatorial description of Lusztig's involution $\eta_\lambda$ on the {\em canonical basis}. Such a description was given by Sch\"utzenberger in type $A$ in terms of tableaux, and  the corresponding procedure is known as {\em evacuation}. The importance of this involution stems from the fact that it exhibits the crystals as self-dual posets, and it corresponds to the action of the longest Weyl group element on an irreducible representation; it also appears in other contexts, such as the recent realization of the category of crystals as a {\em coboundary category} \cite{HK}. Our description of Lusztig's involution is very similar to that of the evacuation map. The main ingredient in defining the latter map, namely Sch\"utzenberger's sliding algorithm (also known as jeu de taquin), is  replaced by Yang-Baxter moves. There is another ingredient, which has to do with  ``reversing'' a $\lambda$-chain and an associated admissible subset, by analogy with reversing the word of a tableau in the definition of the evacuation map. Our construction also leads to a purely combinatorial proof of the fact that the crystals (as defined by our root operators) are self-dual posets. In Section \ref{other-appl}, we present several applications; in particular, we give an intrinsic explanation for the fact that our procedure is an involution.

We will now briefly discuss the relationship between our model and other models for characters. We explained in \cite{LP1} that our model can be viewed as a discrete counterpart to the {\em Littelmann path model}~\cite{Li1,Li2,litcrp,Li3}, which is based on enumerating certain continuous paths in $\hR$. These paths are
constructed recursively, by starting with an initial one, 
and by applying certain root operators. By making specific choices for 
the initial path, one can obtain special cases which have
more explicit descriptions. For instance, a straight line initial path leads to
the {\em Lakshmibai-Seshadri paths} (LS paths) \cite{LS1}; these were introduced before
Littelmann's work, in the context of {\em standard monomial theory}~\cite{LS1}. A model closely related to Littelmann paths is the one due to Gaussent and Littelmann~\cite{GaLi}, which is  based
on {\em LS-galleries}. In \cite{LP,LP1} we discussed in detail the relationship of the alcove path model with Littelmann paths, LS paths, and LS-galleries. We explained the reasons for which the alcove path model is not simply a translation of the Littelmann path model into a different language. We also showed that our model has certain advantages due to its simplicity and combinatorial nature; it also compares favorably in terms of computational complexities (see also Subsection \ref{complexity}).  

The results in this paper highlight new advantages of the alcove path model. For instance, we mentioned above our transparent combinatorial explanation, based on the Yang-Baxter moves, for the independence of the directed colored graph defined by the root operators from the choice of a $\lambda$-chain (Corollary \ref{isom-graphs}). Similarly, it was proved in \cite{Li2} that the directed colored graph structure on Littelmann paths generated by the corresponding root operators is independent of the initial path. However, this proof, which is based on continuous arguments, is less transparent. 

As far as analogs of Sch\"utzenberger's jeu de taquin are concerned, let us first mention the extensions to types $C_n$, $B_n$, and $D_n$ due to Lecouvey and Sheats \cite{lecsc,lecsbd,shesjt}. Let us also note that the only such analog known in the Littelmann path model is the one due to van Leeuwen \cite{leeajt}. The goal of the mentioned paper was to use this analog in order to express in a bijective manner the symmetry of the Littlewood-Richardson rule in the Littelmann path model. 
In a future publication, we will show that van Leeuwen's jeu de taquin realizes precisely the commutator in the category of crystals that was defined in \cite{HK}. 

As far as our combinatorial realization of Lusztig's involution is concerned, let us note that the alcove path model reveals an interesting feature of it, which does not seem to be known even in type $A$. This feature is related to certain Weyl group elements associated with an admissible subset, which we call {\em initial} and {\em final keys} (see Definition \ref{defkeys} and Remark \ref{remkeys}), and which are related to the Demazure character formula in Theorem \ref{demform}. More precisely, Lusztig's involution interchanges the two keys in the sense mentioned in Corollary \ref{conjkeys}. Let us also note that no combinatorial realization of Lusztig's involution is available in the Littelmann path model. However, an explicit description of it is given in \cite{smgrgb} in a different model for characters, which is based on {\em Lusztig's parametrization} and the {\em string parametrization} of the dual canonical basis \cite{baztpm}. Unlike the combinatorial approach in Sch\"utzenberger's evacuation procedure, the involution is now expressed as an affine map whose coefficients are entries of the corresponding Cartan matrix. No intrinsic explanation for the fact that this map is an involution is available. 

We believe that the properties of our model that were investigated in \cite{LP,LP1} as well as in this paper represent just a small fraction of a rich combinatorial structure yet to be explored, which would generalize most of the combinatorics of Young tableaux.  

\section{Preliminaries}\label{prelim}

We recall some background information on finite root systems, affine Weyl groups, Demazure characters, and crystal graphs.

\subsection{Root systems}\label{rootsyst}

Let $G$ be a  connected, simply connected, simple complex Lie group.
Fix a Borel subgroup $B$ and a maximal torus $T$ such that 
$G\supset B\supset T$. As usual, we denote by $B^-$ be the opposite Borel subgroup, while $N$ and $N^-$ are the unipotent radicals of $B$ and $B^-$, respectively.
Let $\mathfrak{g}$, $\h$, $\mathfrak{n}$, and $\mathfrak{n}^-$ be the complex Lie algebras of $G$, $T$, $N$, and $N^-$, respectively.
Let $r$ be the rank of the Cartan subalgebra $\h$.
Let $\Phi\subset \h^*$ be the 
corresponding irreducible {\it root system}, and let $\hR\subset \h^*$ be the real span of the roots.
Let $\Phi^+\subset \Phi$ be the set 
of positive roots corresponding to our choice of $B$. 
Then $\Phi$ is the disjoint union of $\Phi^+$ and $\Phi^- := -\Phi^+$.
We write $\alpha>0$ (respectively, $\alpha<0$) for $\alpha\in\Phi^+$ (respectively, 
$\alpha\in\Phi^-$), and we define ${\rm sgn}(\alpha)$ to be $1$ (respectively $-1$). 
We also use the notation $|\alpha|:={\rm sgn}(\alpha)\alpha$. 
Let $\alpha_1,\dots,\alpha_r\in\Phi^+$ be the corresponding 
{\it simple roots}, which form a basis of $\hR$.
Let $\langle\,\cdot\,,\,\cdot\,\rangle$ denote the nondegenerate scalar product on $\hR$ induced by
the Killing form.  
Given a root $\alpha$, the corresponding {\it coroot\/} is $\alpha^\vee := 2\alpha/\langle\alpha,\alpha\rangle$.  The collection of coroots
$\Phi^\vee:=\{\alpha^\vee \mid \alpha\in\Phi\}$ forms the 
{\it dual root system.}

The {\it Weyl group\/} $W\subset\mathrm{Aut}(\hR)$ 
of the Lie group $G$ is generated by the reflections 
$s_{\alpha}: \hR \to \hR$, for $\alpha\in\Phi$,
given by 
$$
s_\alpha:  \lambda \mapsto \lambda - \langle\lambda,\alpha^\vee\rangle\,\alpha.
$$  
In fact, the Weyl group $W$ is generated by the
{\it simple reflections\/} $s_1,\dots,s_r$ corresponding 
to the simple roots $s_i := s_{\alpha_i}$, subject to the 
{\it Coxeter relations:}
$$(
s_i)^2=1
\quad\textrm{and}\quad 
(s_i s_j)^{m_{ij}}=1\quad\textrm{for any }i,j\in\{1,\dots,r\},
$$
where $m_{ij}$ is half of the order of the dihedral subgroup generated 
by $s_i$ and $s_j$.
An expression of a Weyl group element $w$ as a product 
of generators $w=s_{i_1}\cdots s_{i_l}$
which has minimal length is called a {\it reduced decomposition\/}
for $w$; its length $\ell(w)=l$ is called the {\it length\/} of $w$.
The Weyl group contains a unique {\it longest element\/} $w_\circ$
with maximal length $\ell(w_\circ)=\#\Phi^+$.
For $u,w\in W$, we say that $u$ {\it covers\/} $w$, and write $u\gtrdot w$,
if $w=u s_{\beta}$, for some $\beta\in\Phi^+$, and $\ell(u)=\ell(w)+1$.
The transitive closure ``$>$'' of the relation ``$\gtrdot$'' is called
the {\it Bruhat order\/} on $W$.

The {\it weight lattice\/} $\Lambda$ is given by
\begin{equation}
\Lambda:=\{\lambda\in \hR \mid \langle\lambda,\alpha^\vee\rangle\in\Z
\textrm{ for any } \alpha\in\Phi\}.
\label{eq:weight-lattice}
\end{equation}
The weight lattice $\Lambda$ is generated by the 
{\it fundamental weights\/}
$\omega_1,\dots,\omega_r$, which are defined as the elements of the dual basis to the 
basis of simple coroots, i.e., $\langle\omega_i,\alpha_j^\vee\rangle=\delta_{ij}$.
The set $\Lambda^+$ of {\it dominant weights\/} is given by
$$
\Lambda^+:=\{\lambda\in\Lambda \mid \langle\lambda,\alpha^\vee\rangle\geq 0
\textrm{ for any } \alpha\in\Phi^+\}.
$$

Let $\rho:=\omega_1+\cdots+\omega_r=\frac{1}{2}\sum_{\beta\in\Phi^+}\beta$.
The {\it height\/} of a coroot $\alpha^\vee\in\Phi^\vee$ is 
$\langle\rho,\alpha^\vee\rangle =  c_1+\cdots+c_r$ if 
$\alpha^\vee=c_1\alpha_1^\vee+\cdots + c_r \alpha_r^\vee$.
Since we assumed that $\Phi$ is irreducible, there is 
a unique {\it highest coroot\/} $\theta^\vee\in\Phi^\vee$ that has 
maximal height.  (In other words, $\theta^\vee$ is the highest root
of the dual root system $\Phi^\vee$.  It should not be confused with 
the coroot of the highest root of $\Phi$.)
We will also use the {\it Coxeter number}, that 
can be defined as $h:=\langle\rho,\theta^\vee\rangle+1$.

\subsection{Affine Weyl groups}\label{aff}

In this subsection, we remind a few basic facts about 
affine Weyl groups and alcoves, cf.  
Humphreys~\cite[Chaper~4]{Hum} for more details.

Let $\Waff$ be  the {\it affine Weyl group\/} for the 
Langland's dual group $G^\vee$.
The affine Weyl group $\Waff$ is generated by the affine reflections 
$s_{\alpha,k}: \hR \to \hR$, for $\alpha\in\Phi$ and $k\in\Z$, 
that reflect the space $\hR$ with respect to the affine hyperplanes
\begin{equation}
H_{\alpha,k} := \{\lambda\in \hR \mid \langle\lambda,\alpha^\vee\rangle=k\}.
\label{eq:H-alpha-k}
\end{equation}
Explicitly, the affine reflection $s_{\alpha,k}$ is given by 
$$
s_{\alpha,k}: \lambda \mapsto 
s_{\alpha}(\lambda) + k\,\alpha =
\lambda - (\langle\lambda,\alpha^\vee\rangle-k)\,\alpha.
$$
The hyperplanes $H_{\alpha,k}$ divide the real vector space $\hR$ into open
regions, called {\it alcoves.} 
Each alcove $A$ is given by inequalities of the form
$$
A:=\{\lambda\in \hR \mid m_\alpha<\langle\lambda,\alpha^\vee\rangle<m_\alpha+1
\textrm{ for all } \alpha\in\Phi^+\},
$$
where $m_\alpha=m_\alpha(A)$, $\alpha\in\Phi^+$, are some integers.

A proof of the following important property of the affine Weyl group 
can be found, e.g., in~\cite[Chapter~4]{Hum}.

\begin{lemma}
The affine Weyl group $\Waff$ acts simply transitively 
on the collection of all alcoves.
\label{lem:simply-transitively}
\end{lemma}

The {\it fundamental alcove\/} $A_\circ$ is given by 
$$
A_\circ :=\{\lambda\in \hR \mid 0<\langle\lambda,\alpha^\vee\rangle<1 \textrm{ for all }
\alpha\in\Phi^+\}.
$$
Lemma~\ref{lem:simply-transitively} implies that, for any alcove $A$, 
there exists a unique element $v_A$ of the affine Weyl group $\Waff$
such that $v_A(A_\circ) = A$.  Hence the map $A\mapsto v_A$ is a one-to-one
correspondence between alcoves and elements of the affine Weyl group.

Recall that $\theta^\vee\in\Phi^\vee$ is the highest coroot. 
Let $\theta\in\Phi^+$ be the corresponding root,
and let $\alpha_0:=-\theta$.
The fundamental alcove $A_\circ$ is, in fact, the simplex given by
\begin{equation}
A_\circ =\{\lambda\in \hR \mid 0<\langle\lambda,\alpha_i^\vee\rangle \textrm{ for }
i=1,\dots,r, \textrm{ and }\langle\lambda,\theta^\vee\rangle<1\},
\label{eq:fund-alcove}
\end{equation}
Lemma~\ref{lem:simply-transitively} also implies that the affine Weyl group 
is generated by the set of reflections $s_0,s_1,\dots,s_r$ with respect 
to the walls of the fundamental alcove $A_\circ$, where 
$s_0 := s_{\alpha_0,-1}$ and $s_1,\dots,s_r\in W$ 
are the simple reflections $s_i=s_{\alpha_i,0}$. Like the Weyl group, the affine Weyl group $\Waff$ is a Coxeter group. As in the case of the Weyl group, a decomposition $v=s_{i_1}\cdots s_{i_l}\in \Waff$ is called {\it reduced\/} if it has minimal length;
its length $\ell(v)=l$ is called the length of $v$.

We say that two alcoves $A$ and $B$ are {\it adjacent} 
if $B$ is obtained by an affine reflection of $A$ with respect to one of its 
walls.  In other words, two alcoves are adjacent if they are
distinct and have a common wall.  
For a pair of adjacent alcoves, let us write 
$A\stackrel{\beta}\longrightarrow B$ if the common wall of $A$ and $B$ 
is of the form $H_{\beta,k}$ and the root $\beta\in\Phi$ points 
in the direction from $A$ to $B$.  

 Let $Z$ be the set of the elements of the lattice $\Lambda/\hvee$ 
that do not belong to any affine hyperplane $H_{\alpha,k}$ (recall that $h$ is the Coxeter number). Each alcove $A$ contains precisely one element $\zeta_A$ of the set $Z$ (cf.\ \cite{Kost,LP}); this will be called the {\em central point} of $A$. In particular, $\zeta_{A_\circ}=\rho/\hvee$. 

\begin{proposition}\cite{LP}\label{adjalcoves} For a pair of adjacent alcoves $A\stackrel{\alpha}\longrightarrow B$, we have $\zeta_B -\zeta_A = \alpha/\hvee$.
\end{proposition}

\subsection{Demazure characters}

The {\it generalized flag variety\/} $G/B$ is a smooth projective
variety.  It decomposes into a disjoint union of {\it Schubert cells\/} 
$X_w^\circ := BwB/B$ indexed by elements $w\in W$ of the Weyl group.
The closures of Schubert cells $X_w := \overline{X_w^\circ}$
are called {\it Schubert varieties.} 
We have $u>w$ in the Bruhat order (defined above) 
if and only if $X_u\supset X_w$. 
Let $\O_{X_w}$ be the structure sheaf of the Schubert variety $X_w$. Let $\mathcal{L}_\lambda$ be the line bundle over $G/B$ associated with the
weight $\lambda$, that is, $\mathcal{L}_\lambda:= G\times_B \C_{-\lambda}$, where $B$ acts on $G$ by right multiplication, and the
$B$-action on $\C_{-\lambda}=\C$ corresponds to the character determined by $-\lambda$. 
(This character of $T$ extends to $B$ by defining it to be 
identically one on the commutator subgroup $[B,B]$.) 

For a dominant weight $\lambda\in \Lambda^+$,
let $V_\lambda$ denote the finite dimensional irreducible representation of 
the Lie group $G$ with highest weight $\lambda$. 
For $\lambda\in\Lambda^+$ and $w\in W$, 
the {\it Demazure module\/} $V_{\lambda,w}$ is the $B$-module that 
is dual to the space of global 
sections of the line bundle $\L_\lambda$ on the Schubert variety $X_w$:
\begin{equation}\label{demaz}
V_{\lambda,w} := H^0(X_w,\L_\lambda)^*.  
\end{equation}
For the longest Weyl group element $w=w_\circ$,
the space $V_{\lambda,w_\circ} = H^0(G/B,\L_\lambda)^*$ has the structure of
a $G$-module.  The classical {\it Borel-Weil theorem\/} says that 
$V_{\lambda,w_\circ}$ is isomorphic to the irreducible $G$-module $V_\lambda$.

Let $\Z[\Lambda]$ be the group algebra of the weight lattice $\Lambda$, which is isomorphic to the representation ring of $T$.  The algebra $\Z[\Lambda]$ has
a $\Z$-basis of formal exponents $\{e^\lambda \mid \lambda\in\Lambda\}$ with
multiplication $e^\lambda\cdot e^\mu := e^{\lambda+\mu}$; in other words,
$\Z[\Lambda]=\Z[e^{\pm\omega_1},\cdots,e^{\pm\omega_r}]$ is the algebra of
Laurent polynomials in  $r$ variables.  
The formal characters of the modules $V_{\lambda,w}$, called {\it Demazure characters}, 
are given by $ch(V_{\lambda,w})=\sum_{\mu\in\Lambda} 
m_{\lambda,w}(\mu)\,e^\mu\in\Z[\Lambda]$,
where $m_{\lambda,w}(\mu)$ is the multiplicity of the weight $\mu$ in 
$V_{\lambda,w}$.
These characters generalize the characters of the irreducible representations
$ch(V_\lambda)=ch(V_{\lambda,w_\circ})$.
Demazure~\cite{Dem} gave a formula expressing the characters $ch(V_{\lambda,w})$ in terms of certain operators known as {\em Demazure operators}. 

\subsection{Crystal graphs and Lusztig's involution} \label{cryst-inv} Let $U({\mathfrak g})$ be the universal enveloping algebra of the Lie algebra ${\mathfrak g}$. Let $\mathcal{B}$ be the {\em canonical basis} of $U({\mathfrak n}^-)$, and let ${\mathcal B}_\lambda:={\mathcal B}\cap V_\lambda$ be the canonical basis of the irreducible representation $V_\lambda$ with highest weight $\lambda$. Let $v_\lambda$ and $v_\lambda^{low}$ be the highest and lowest weight vectors in ${\mathcal B}_\lambda$, respectively. Let $\widetilde{E}_i$, $\widetilde{F}_i$, for $i=1,\ldots,r$, be Kashiwara's operators \cite{kascqa,lusiqg}; these are also known as raising and lowering operators, respectively.  The {\em crystal graph} of $V_\lambda$ is the directed colored graph on ${\mathcal B}_\lambda$ defined by  arrows $x\rightarrow y$ colored $i$ for each $\widetilde{F}_i(x)=cy+\,{\rm lower}\;\,{\rm terms}$, or, equivalently, for each $\widetilde{E}_i(y)=cx+\,{\rm lower}\;\,{\rm terms}$, with $c$ a constant. (In fact, Kashiwara introduced the notion of a crystal graph of an $U_q({\mathfrak g})$-representation, where $U_q({\mathfrak g})$ is the Drinfeld-Jimbo $q$-deformation of $U({\mathfrak g})$, also known as a {\em quantum group}; using the quantum deformation, one can associate a crystal graph to a ${\mathfrak g}$-representation.) One can also define partial orders $\preceq_i$ on ${\mathcal B}_\lambda$ by\[ x\preceq_iy\;\;\;\mathrm{if}\;\;\;x=\widetilde{F}_i^k(y)\;\:\textrm{for some }\,k\ge 0\,.\]
We let $\preceq$ denote the partial order generated by all partial orders $\preceq_i$, for $i=1,\ldots,r$. The poset $({\mathcal B}_\lambda,\preceq)$ has maximum $v_\lambda$ and minimum $v_\lambda^{low}$. 

In order to proceed, we need the following general setup. Let $V$ be a module over an associative algebra $U$ and $\sigma$ an automorphism of $U$. The twisted $U$-module $V^\sigma$ is the same vector space $V$ but with the new action $u*v:=\sigma(u) v$ for $u\in U$ and $v\in V$. Clearly, $V^{\sigma\tau}=(V^\sigma)^\tau$ for every two automorphisms $\sigma$ and $\tau$ of $U$. Furthermore, if $V$ is a simple $U$-module, then so is $V^\sigma$. In particular, if $U=U({\mathfrak g})$ and $V=V_\lambda$, then $(V_\lambda)^\sigma$ is isomorphic to $V_{\sigma(\lambda)}$ for some dominant weight $\sigma(\lambda)$. Thus there is an isomorphism of vector spaces $\sigma_\lambda\::\:V_\lambda\rightarrow V_{\sigma(\lambda)}$ such that 
\[\sigma_\lambda(uv)=\sigma(u)\sigma_\lambda(v)\,,\;\;\;\;u\in U({\mathfrak g})\,,\;v\in V_\lambda\,.\]
By Schur's lemma, $\sigma_\lambda$ is unique up to a scalar multiple. 

The longest Weyl group element $w_\circ$ defines an involution on the simple roots by $\alpha_i\mapsto \alpha_{i^*}:=-w_\circ(\alpha_i)$. Consider the automorphisms of $U({\mathfrak g})$ defined by
\begin{align}
&\phi(E_i)=F_i\,,\;\;\;\;\;\;\phi(F_i)=E_i\,,\;\;\;\;\;\phi(H_i)=-H_i\,,\label{defphi}\\
&\psi(E_i)=E_{i^*}\,,\;\;\;\;\psi(F_i)=F_{i^*}\,,\;\;\;\;\psi(H_i)=H_{i^*}\,,\label{defpsi}
\end{align}
and $\eta:=\phi\psi$. Clearly, these three automorphisms together with the identity automorphism form a group isomorphic to $\Z/2\Z\times\Z/2\Z$. It also easily follows from (\ref{defphi})-(\ref{defpsi}) that 
\[\phi(\lambda)=\psi(\lambda)=-w_\circ(\lambda)\,,\;\;\;\;\eta(\lambda)=\lambda\,.\]
We can normalize each of the maps $\phi_\lambda$, $\psi_\lambda$, and $\eta_\lambda$ by the requirement that
\begin{equation}\label{phipsieta0}\phi_\lambda(v_\lambda)=v_{-w_\circ(\lambda)}^{low}\,,\;\;\;\;\psi_\lambda(v_\lambda)=v_{-w_\circ(\lambda)}\,,\;\;\;\;\eta_\lambda(v_\lambda)=v_\lambda^{low}\,.\end{equation}
(Of course, we also set $\mathrm{Id}_\lambda$ to be the identity map on $V_\lambda$.) By \cite[Proposition 21.1.2]{lusiqg}, cf. also \cite[Proposition 7.1]{bazcbq}, we have the following result.

\begin{proposition}\cite{bazcbq,lusiqg}\label{phipsieta} {\rm (1)} Each of the maps $\phi_\lambda$ and $\psi_\lambda$ sends ${\mathcal B}_\lambda$ to ${\mathcal B}_{-w_\circ(\lambda)}$, while $\eta_\lambda$ sends ${\mathcal B}_\lambda$ to itself.

{\rm (2)} For every two (not necessarily distinct) elements $\sigma$, $\tau$ of the group $\{\mathrm{Id},\phi,\psi,\eta\}$, we have $(\sigma\tau)_\lambda=\sigma_{\tau(\lambda)}\tau_\lambda$. In particular, the map $\eta_\lambda$ is an involution.

{\rm (3)} For every $i=1,\ldots,r$, we have
\begin{equation}\label{phipsieta3}\phi_\lambda\widetilde{F}_i=\widetilde{E}_i\phi_\lambda\,,\;\;\;\;\psi_\lambda\widetilde{F}_i=\widetilde{F}_{i^*}\psi_\lambda\,,\;\;\;\;\eta_\lambda\widetilde{F}_i=\widetilde{E}_{i^*}\eta_\lambda\,.\end{equation}
In particular, the poset $({\mathcal B}_\lambda,\preceq)$ is self-dual, and $\eta_\lambda$ is the corresponding antiautomorphism. 
\end{proposition}

Berenstein and Zelevinsky \cite{bazcbq} showed that, in type $A_{n-1}$ (that is, in the case of the Lie algebra ${\mathfrak s}{\mathfrak l}_{n}$), the operator $\eta_\lambda$ is given by Sch\"utzenberger's {\em evacuation} procedure for semistandard Young tableaux (see e.g. \cite{fulyt}). More precisely, it is known that, for each partition $\lambda=(\lambda_1\ge\lambda_2\ge\ldots\ge\lambda_{n-1}\ge 0)$, the semistandard Young tableaux of shape $\lambda$ and entries $1,\ldots,n$ parametrize the  the canonical basis ${\mathcal B}_\lambda$ of $V_\lambda$. Hence, we can transfer the action of $\eta_\lambda$ on ${\mathcal B}_\lambda$ to an action on the corresponding tableaux. As mentioned above, the latter action coincides with Sch\"utzenberger's evacuation map. One way to realize this map on a tableau $T$ is the following three-step procedure.
\begin{enumerate}
\item Rotate the tableau $180^\circ$, such that its row/column words get reversed. 
\item Complement the entries via the map $i\mapsto w_\circ(i)=n+1-i$, where $w_\circ$ is the longest element in the symmetric group $S_{n}$.
\item Apply {\em jeu de taquin} to construct the {\em rectification} of the skew tableau obtained in the previous step, that is, successively apply Sch\"utzenberger's {\em sliding algorithm} for the inside corners of the mentioned tableau.
\end{enumerate}
For convenience, we will call these steps: REVERSE, COMPLEMENT, SLIDE. They are illustrated in Figure \ref{fig:evac} below.
\begin{figure}[ht]
\mbox{\epsfig{file=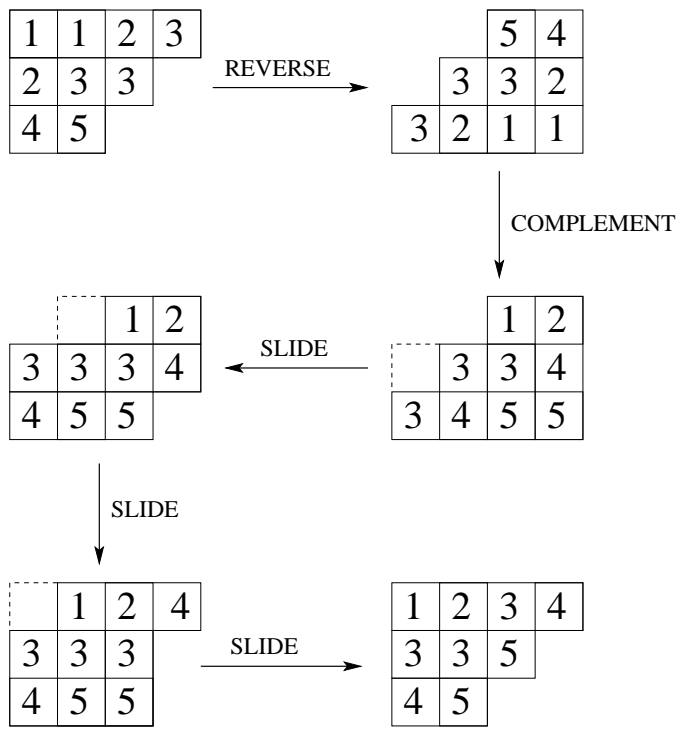}}
\caption{The evacuation map.}
\label{fig:evac}
\end{figure}

\section{The Alcove Path Model}\label{alc-path}

In this section, we recall the model for the irreducible characters of semisimple Lie algebras that we introduced in \cite{LP,LP1}. We refer to these papers for more details, including the proofs of the results mentioned below. Although some of these results hold for infinite root systems (cf. \cite{LP1}), the setup in this paper is that of a finite irreducible root system, as discussed in Section \ref{prelim}. 

Our model is conveniently phrased in terms of several sequences, so let us mention some related notation. Given a totally ordered index set $I=\{i_1<i_2<\ldots<i_n\}$, a sequence $(a_{i_1},a_{i_2},\ldots,a_{i_n})$ is sometimes abbreviated to $\{a_j\}_{j\in I}$. We also let $[n]:=\{1,\,2,\,\ldots,\,n\}$. 

\subsection{$\lambda$-chains}\label{lamch} The affine translations by weights preserve the set
of affine hyperplanes $H_{\alpha,k}$, 
cf.~(\ref{eq:weight-lattice}) and~(\ref{eq:H-alpha-k}).
It follows that these affine translations map alcoves to alcoves.
Let $A_\lambda=A_\circ + \lambda$
be the alcove obtained by the affine translation of the fundamental alcove
$A_\circ$ by a weight $\lambda\in\Lambda$.  Let $v_\lambda$
be the corresponding element of $\Waff$, i.e,.
$v_{\lambda}$ is defined by $v_\lambda(A_\circ) = A_\lambda$.
Note that the element $v_\lambda$ may not be an affine translation
itself.

Let us now fix a dominant weight $\lambda$. Let $v\mapsto \bar v$ be the homomorphism $\Waff\to W$ defined by ignoring the affine translation. In other words, $\bar s_{\alpha,k} = s_\alpha\in W$.

\begin{definition}
A {\it $\lambda$-chain of roots\/} is a sequence of positive roots $(\beta_1,\dots,\beta_n)$ which is determined as indicated below by a reduced decomposition $v_{-\lambda} = s_{i_1}\cdots s_{i_n}$ 
of $v_{-\lambda}$ as a product of generators of $\Waff$:
$$
\beta_1 = \alpha_{i_1}, \  \beta_2=\bar s_{i_1}(\alpha_{i_2}), \ 
\beta_3 = \bar s_{i_1} \bar s_{i_2} (\alpha_{i_3}), \dots, \ 
\beta_n= \bar s_{i_1} \cdots \bar s_{i_{n-1}} (\alpha_{i_n})\,.
$$
When the context allows, we will abbreviate ``$\lambda$-chain of roots'' to ``$\lambda$-chain''.
The {\it $\lambda$-chain of reflections\/} associated 
with the above  $\lambda$-chain of roots is the sequence 
$(\widehat{r}_1,\dots,\widehat{r}_n)$ of affine reflections in $\Waff$ given by
$$
\widehat{r}_1 = s_{i_1}, \ \widehat{r}_2 = s_{i_1} s_{i_2} s_{i_1}, \ 
\widehat{r}_3 =  s_{i_1} s_{i_2} s_{i_3} s_{i_2} s_{i_1}, \ \dots, \
\widehat{r}_n =  s_{i_1} \cdots s_{i_n} \cdots s_{i_1}. 
$$
\label{def:lambda-chain}
\end{definition}

We will present two equivalent definitions of a $\lambda$-chain of roots.

\begin{definition}
An {\it alcove path\/} is a sequence of alcoves
$(A_0,A_1,\dots,A_n)$ such that $A_{i-1}$ and $A_i$ are adjacent, for
$i=1,\dots,n$.
We say that an alcove path is {\it reduced\/} if it has minimal 
length among all alcove paths from $A_0$ to $A_n$.
\end{definition}

Given a finite sequence of roots $\Gamma=(\beta_1,\ldots,\beta_n)$, we define the sequence of integers $(l_1^\emptyset,\ldots,l_n^\emptyset)$ by $l_i^\emptyset:= \#\{j<i\mid \beta_j = \beta_i\}$, for $i=1,\dots,n$. We also need the following two conditions on $\Gamma$.
\begin{enumerate}
\item[(R1)] The number of occurrences of any positive root $\alpha$ in $\Gamma$ is 
$\langle \lambda,\alpha^\vee \rangle$. 
\item[(R2)] For each triple of positive roots $(\alpha,\,\beta,\,\gamma)$ with $\gamma^\vee=\alpha^\vee+\beta^\vee$, the subsequence of $\Gamma$ consisting of $\alpha,\,\beta,\,\gamma$ is a concatenation of pairs $(\alpha,\gamma)$ and $(\beta,\gamma)$ (in any order). 
\end{enumerate}

\begin{theorem}\cite{LP}\label{equivdef}
The following statements are equivalent.
\begin{enumerate}
\item[(a)] The sequence of roots $\Gamma=(\beta_1,\ldots,\beta_n)$ is a $\lambda$-chain, and $(\widehat{r}_1,\ldots,\widehat{r}_n)$ is the associated $\lambda$-chain of reflections.
\item[(b)] We have a reduced alcove path $A_0\stackrel{-\beta_1}\longrightarrow \cdots
\stackrel{-\beta_n}\longrightarrow A_{n}$ from $A_0=A_\circ$ to $A_n=A_{-\lambda}$, and $\widehat{r}_i$ is the affine reflection in the common wall of $A_{i-1}$ and $A_i$, for $i=1,\dots,n$.
\item[(c)] The sequence $\Gamma$ satisfies conditions {\rm (R1)} and {\rm (R2)} above, and $\widehat{r}_i=s_{\beta_i,-l_i^\emptyset}$, for $i=1,\dots,n$.
\end{enumerate}
\end{theorem}

We now describe a particular choice of a $\lambda$-chain. First note that constructing a $\lambda$-chain amounts to defining a total order on the index set 
\[
I:=\{(\alpha,k)\mid  \alpha\in\Phi^+,\,0\le k<\langle \lambda,\alpha^\vee \rangle\}\,,\]
such that condition (R2) above holds, where the sequence  $\Gamma=\{\beta_i\}_{i\in I}$ is defined by $\beta_i=\alpha$ for $i=(\alpha,k)$.  Fix a total
order on the set of simple roots $\alpha_1<\alpha_2<\ldots<\alpha_r$. For each
$i=(\alpha,k)$ in $I$, let
$\alpha^\vee=c_1\alpha_1^\vee+\ldots+c_r\alpha_r^\vee$, and define the vector 
\[
v_i:=\frac{1}{\langle \lambda,\alpha^\vee \rangle}(k,c_1,\ldots,c_r)
\]
in $\Q^{r+1}$. 
It turns out that the map $i\mapsto v_i$ is injective. Hence, we can define a total order on $I$ by $i<j$ iff $v_i<v_j$ in the lexicographic order on $\Q^{r+1}$. 

\begin{proposition}\cite{LP1}\label{constr-lambdachain} Given the total order on $I$ defined above, the sequence $\{\beta_i\}_{i\in I}$ defined by $\beta_i=\alpha$ for $i=(\alpha,k)$ is a $\lambda$-chain.
\end{proposition}

\subsection{Admissible subsets} For the remainder of this section, we fix a $\lambda$-chain $\Gamma=(\beta_1,\ldots,\beta_n)$. Let $r_i:=s_{\beta_i}$. We now define the centerpiece of our combinatorial model for characters, which is our generalization of semistandard Young tableaux in type $A$. 

\begin{definition}
An {\em admissible subset} is a subset of $[n]$ (possibly empty), that is, $J=\{j_1<j_2<\ldots<j_s\}$, such that we have the following saturated chain in the Bruhat order on $W$:
\[1\lessdot r_{j_1}\lessdot r_{j_1}r_{j_2}\lessdot \ldots \lessdot r_{j_1}r_{j_2}\ldots r_{j_s}\,.\] We denote by ${\mathcal A}(\Gamma)$ the collection of all admissible subsets corresponding to our fixed $\lambda$-chain $\Gamma$. Given an admissible subset $J$, we use the notation
\[\mu(J):=-\widehat{r}_{j_1}\ldots\widehat{r}_{j_s}(-\lambda)\,,\;\;\;\;\;w(J):=r_{j_1}\ldots r_{j_s}\,.\]
We call $\mu(J)$ the {\em weight} of the admissible subset $J$.
\end{definition}

\begin{theorem}\cite{LP,LP1}\label{charform} {\rm (1)} We have the following character formula:
\[ch(V_\lambda)=\sum_{J\in {\mathcal A}(\Gamma)}e^{\mu(J)}\,.\]

{\rm (2)} More generally, the following Demazure character formula holds for any $u\in W$:
$$
ch(V_{\lambda,u}) =
\sum_J e^{-u\, \widehat{r}_{j_1}\cdots \widehat{r}_{j_s}(-\lambda)}\,,
$$
where the summation is over all subsets $J=\{j_1<\cdots<j_s\}\subseteq[n]$ such that
$$
u \gtrdot u\,  r_{j_1} \gtrdot 
u \, r_{j_1}  r_{j_2} \gtrdot \cdots \gtrdot 
u \, r_{j_1}  r_{j_2} \cdots r_{j_s} 
$$
is a saturated decreasing chain in the Bruhat order on the Weyl group $W$.
\end{theorem}

In addition to the above character formulas, a Littlewood-Richardson rule for decomposing tensor products of irreducible representations is also presented in terms of our model in \cite{LP1}. 

\begin{example}\label{simpleex}  Consider the Lie algebra $\mathfrak{sl}_3$ of type $A_2$. The corresponding root system $\Phi$ can be realized inside the vector space $V:=\R^3/\R(1,1,1)$ as $\Phi=\{\alpha_{ij}:=\varepsilon_i-\varepsilon_j \mid i\ne j,\,1\le i,j\le 3\}$, where $\varepsilon_1,\varepsilon_2,\varepsilon_3\in V$ are the images of the coordinate vectors in $\R^3$. The reflection $s_{\alpha_{ij}}$ is denoted by $s_{ij}$. The simple roots are $\alpha_{12}$ and $\alpha_{23}$, while $\alpha_{13}=\alpha_{12}+\alpha_{23}$ is the other positive root.
Let $\lambda=\omega_1=\varepsilon_1$ be the first fundamental weight.  
In this case, there is only one $\lambda$-chain
$(\beta_1,\beta_2)=(\alpha_{12},\alpha_{13})$. There are 3 admissible subsets: $\emptyset, \{1\}, \{1,2\}$.
The subset $\{2\}$ is not admissible because the reflection 
$s_{13}$ does not cover the identity element.
We have $(l_1^\emptyset,l_2^\emptyset) = (0,0)$.
Theorem~\ref{charform} gives the following expression
for the character of $V_{\omega_1}$:
$$
ch(V_{\omega_1}) = e^{\omega_1} + e^{s_{12}(\omega_1)} + 
e^{s_{12} s_{13}(\omega_1)}.
$$
\end{example}

In Subsections \ref{subs:galleries} and \ref{chainroots}, we present two alternative ways of viewing admissible subsets, which are closely related to the equivalent definitions of $\lambda$-chains in Theorem \ref{equivdef} (b) and (c).

\subsection{Computational complexities}\label{complexity}

In this subsection, we compare the computational complexity of our
model with that of LS-paths constructed via root operators.

Fix a root system of rank $r$ with $N$ positive roots, a dominant weight
$\lambda$, and a Weyl group element $u$ of length $l$. We want to determine the
character of the Demazure module $V_{\lambda,u}$. Let $d$ be its dimension, and
let $L$ be the length of the affine Weyl group element $v_{-\lambda}$ (that is,
the number of affine hyperplanes separating the fundamental alcove $A_\circ$
and $A_\circ-\lambda$). Note that $L=2(\lambda, \rho^\vee)$, where
$\rho^\vee=\frac{1}{2}\sum_{\beta\in\Phi^+}\beta^\vee$. We claim that the
complexity of the character formula in Theorem \ref{charform} (2) is $O(d\, l L)$. Indeed, we start by
determining an alcove path via the method underlying Proposition \ref{constr-lambdachain}, which involves
sorting a sequence of $L$ rational numbers. The complexity is $O(L\,\log L)$,
and note that $\log L$ is, in general, much smaller than $d$ (see below for
some examples). Whenever we examine some subword of the word of length $L$ we
fixed at the beginning, we have to check at most $L-1$ ways to add an extra
reflection at the end. On the other hand, in each case, we have to check
whether, upon multiplying by the corresponding nonaffine reflection, the length
decreases by precisely 1. The complexity of the latter operation is $O(l)$,
based on the Strong Exchange Condition \cite[Theorem~5.8]{Hum}. Then, for each
``good'' subword, we have to do a calculation, namely applying at most $2l$
affine reflections to $-\lambda$. In fact, it is fairly easy to implement this
algorithm.

Now let us examine at the complexity of the algorithm 
based on root operators for constructing the LS-paths associated with $\lambda$.
In other words, we
are looking at the complexity of constructing the corresponding crystal graph.
We have to generate the whole crystal graph first, and then figure out which
paths give weights for the Demazure module. For each path, we can apply $r$
root operators. Each path has at most $N$ linear steps, so applying a root
operator has complexity $O(N)$. But now we have to check whether the result is
a path already determined, so we have to compare the obtained path with the
other paths (that were already determined) of the same rank in the crystal
graph (viewed as a ranked poset). This has complexity $O(N M)$, where $M$ is
the maximum number of elements of the same rank. Since we have at most $N+1$
ranks, $M$ is at least $d/(N+1)$. In conclusion, the complexity is $O(d r N
M)$, which is at least $O(d^2 r)$.

Let us get a better picture of how the two results compare. Assume we are in a
classical type, and let us first take $\lambda$ to be the $i$-th fundamental
weight, with $i$ fixed, plus $u=w_\circ$. Clearly $l$ is $O(r^2)$, $L$ is
$O(r)$, and $d$ is $O(r^i)$, so the complexity of our formula is $O(r^{i+3})$.
For LS-paths, we get at least $O(r^{2i+1})$. So the ratio between the
complexity in the model based on LS-paths and our model is at least
$O(r^{i-2})$.

Let us also take $\lambda=\rho$. In this case $d=2^N$, and a simple calculation
shows that $L$ is $O(r^3)$. Our formula has complexity $O(2^N r^5)$, while the
model based on LS-paths has complexity at least $O(2^{2N} r)$. So the ratio
between the complexities is at least $O(2^N/r^4)$, where $N$ is $r(r+1)/2$,
$r^2$, and $r^2-r$ in types $A$, $B/C$, and $D$, respectively.

\subsection{Galleries} \label{subs:galleries}

\begin{definition} A {\it gallery\/} is a sequence 
$\gamma=(F_0=\{0\},A_0=A_\circ,F_1,A_1, F_2, \dots , F_n, A_n, F_{\infty}=\{\mu\})$ 
such that $A_0,\dots,A_n$ are alcoves;
$F_i$ is a codimension one common face of the alcoves $A_{i-1}$ and $A_i$,
for $i=1,\dots,n$; and 
$F_{\infty}$ is a vertex of the last alcove $A_n$. 
The weight $\mu$ is called the {\it weight\/} of the gallery and is denoted by $\mu(\gamma)$. 
The folding operator $\phi_i$ is the operator which acts on a gallery by  leaving its initial segment from $A_0$ to $A_{i-1}$ intact and by reflecting the remaining tail in the affine hyperplane containing the face $F_i$. In other words, we define
$$\phi_i(\gamma):=(F_0,A_0, F_1, A_1, \dots, A_{i-1}, F_i'=F_i,  A_{i}', F_{i+1}', A_{i+1}', \dots,  A_n', F_{\infty}')\,;$$
here $A_j' := \widehat{t}_i(A_j)$ for $j\in\{i,\dots,n\}$, $F_j':=\widehat{t}_i(F_j)$ for $j\in\{i,\dots,n\}\cup\{\infty\}$, and $\widehat{t}_i$ is the affine reflection in the hyperplane containing $F_i$, as in Theorem \ref{equivdef}.
\end{definition}

 The galleries defined above are special cases of the generalized galleries in~\cite{GaLi}.

Recall that our fixed $\lambda$-chain $\Gamma=(\beta_1,\ldots,\beta_n)$ determines a reduced alcove path $A_0=A_\circ\stackrel{-\beta_1}\longrightarrow \cdots
\stackrel{-\beta_n}\longrightarrow A_{n}=A_{-\lambda}$. This alcove path determines, in turn, an obvious gallery 
\[\gamma(\emptyset)=(F_0,A_0,F_1,\dots,F_n, A_n, F_{\infty})\]
 of weight $-\lambda$. 

\begin{definition} 
Given a subset $J= \{ j_1<\cdots< j_s\}\subseteq [n]$, we associate with it the gallery $\gamma(J):=\phi_{j_1}\cdots \phi_{j_s} (\gamma(\emptyset))$. If $J$ is an admissible subset, we call $\gamma(J)$ an {\it admissible gallery}. 
\end{definition}

\begin{remarks} (1) The weight of the gallery $\gamma(J)$, i.e. $\mu(\gamma(J))$, is $-\mu(J)$.

(2) In order to define the gallery $\gamma(J)$, we augmented the index set $[n]$ corresponding to the fixed $\lambda$-chain by adding a new minimum $0$ and a new maximum $\infty$. The same procedure is applied when the initial index set is an arbitrary (finite) totally ordered set.
\end{remarks}

\subsection{Chains of roots} \label{chainroots}

\begin{definition} A {\em chain of roots} is an object of the form 
\begin{equation}\label{defgamma}\Gamma=((\gamma_1,\gamma_1'),\ldots,(\gamma_n,\gamma_n'),\gamma_\infty) \,,\end{equation}
where $(\gamma_i,\gamma_i')$ are pairs of roots with $\gamma_i'=\pm\gamma_i$, for $i=1,\dots,n$, and $\gamma_\infty$ is a weight. Given a chain of roots $\Gamma$ and $i$ in $[n]$, we let $t_i:=s_{\gamma_i}$ and we define 
\[\phi_i(\Gamma):=((\delta_1,\delta_1'),\ldots,(\delta_n,\delta_n'),\delta_\infty)) \,,\]
where $\delta_\infty:={t}_i(\gamma_\infty)$ and 
\[(\delta_j,\delta_j'):=\casethree{(\gamma_j,\gamma_j')}{j<i}{(\gamma_j,t_i(\gamma_j'))}{j=i}{(t_i(\gamma_j),t_i(\gamma_j'))}{j>i}\]
\end{definition}

Our fixed $\lambda$-chain $\Gamma=(\beta_1,\ldots,\beta_n)$ determines the chain of roots
\[\Gamma(\emptyset):=((\beta_1,\beta_1), \ldots, (\beta_n,\beta_n),\rho) \,;\]
recall that  $\rho$ was defined in Subsection \ref{rootsyst}.
 
\begin{definition} 
Given a subset $J= \{ j_1<\cdots< j_s\}\subseteq [n]$, we associate with it the chain of roots $\Gamma(J):=\phi_{j_1}\cdots \phi_{j_s} (\Gamma(\emptyset))$. If $J$ is an admissible subset, we call $\Gamma(J)$ an {\it admissible folding} (of $\Gamma(\emptyset)$). 
\end{definition}

\begin{remark} We can also define folding operators on subsets $J$ of $[n]$ by $\phi_i\::\:J\mapsto J\triangle\{i\}$, where $\triangle$ denotes the symmetric difference of sets. The folding operators $\phi_i$ on $J$, $\gamma(J)$, and $\Gamma(J)$ are commuting involutions (for $J\subseteq[n]$), and their actions are compatible. Throughout this paper, we use $J$, $\gamma(J)$, and $\Gamma(J)$ interchangeably. We will call the elements of $J$ the {\em folding positions} in $\gamma(J)$ and $\Gamma(J)$. \end{remark}

Given a fixed subset $J$ of $[n]$, we will now discuss the relationship between the gallery $\gamma(J)$ and the chain of roots $\Gamma(J)$. 

Let $\gamma=(F_0,A_0,F_1,\dots,F_n, A_n, F_{\infty})$ be an arbitrary gallery. Let $\widehat{t}_i$ be the affine reflection in the common wall of $A_{i-1}$ and $A_i$, as usual. We associate with $\gamma$ a chain of roots $\Gamma(\gamma)=((\gamma_1,\gamma_1'),\ldots,(\gamma_n,\gamma_n'),\gamma_\infty)$ as follows:
 \begin{equation}\label{gal-fol1}\gamma_i:=h(\zeta_{A_{i-1}}-\zeta_{\widehat{t}_i(A_{i-1})})\,,\;\;\;\;\gamma_i':=h(\zeta_{\widehat{t}_i(A_{i})}-\zeta_{A_{i}})\,,\;\;\;\;\gamma_\infty:=h(\zeta_{A_n}-\mu(\gamma))\,;\end{equation}
here $h$ is the Coxeter number, $i=1,\ldots,n$, and $\zeta_A$ is the central point of the alcove $A$, as defined in Subsection \ref{aff}. By Proposition \ref{adjalcoves}, we have
 \begin{equation}\label{gal-fol2}\widehat{t}_i(A_{i-1}) \stackrel{\gamma_i}\longrightarrow A_{i-1}\,,\;\;\;\;\;A_i  \stackrel{\gamma_i'}\longrightarrow\widehat{t}_i(A_i)\,.\end{equation}
On the one hand, $\Gamma(\gamma)$ uniquely determines the gallery $\gamma$. On the other hand,  we have $\Gamma(J)=\Gamma(\gamma(J))$. 

\begin{remark}
In \cite{LP1}, we also associated with an admissible subset $J$ a certain piecewise-linear path. This is closely related to $\gamma(J)$ and $\Gamma(J)$; essentially, it is obtained from the path joining the central points of the alcoves in the gallery $\gamma(\emptyset)$ via the folding operators used to construct $\gamma(J)$ from $\gamma(\emptyset)$. However, this path is {\em not} a Littelmann path in general.
\end{remark}

\subsection{Combinatorial properties}\label{comb-adm-seq} Let $J$ be a fixed admissible subset, and let 
\[\gamma(J)=(F_0,A_0,F_1,\dots,F_n, A_n, F_{\infty})\,,\;\;\;\;\Gamma(J)=((\gamma_1,\gamma_1'),\ldots,(\gamma_n,\gamma_n'),\gamma_\infty)\,.\]
Let us also fix a simple root $\alpha_p$. We associate with $J$ the sequence of integers $L(J)=(l_1,\ldots,l_n)$ defined by $F_i\subset H_{-|\gamma_i|,l_i}$ for $i=1,\ldots,n$. Note that $L(\emptyset)=(l_1^\emptyset,\ldots,l_n^\emptyset)$, as defined in Subsection \ref{lamch}. We also define $l_\infty^p:=\langle\mu(J),\alpha_p^\vee\rangle$, which means that $F_\infty\subset H_{-\alpha_p,l_\infty^p}$. Finally, we let
\begin{equation}\label{ilmp}I(J,p):=\{i\in[n]\mid \gamma_i=\pm\alpha_p\}\,,\;\;\;\;L(J,p):=(\{l_i\}_{i\in I(J,p)},\,l_\infty^p)\,,\;\;\;\; M(J,p):=\max\,{L}(J,p)\,.\end{equation}
 It turns out that $M(J,p)\ge 0$. 

Let $I(J,p)=\{i_1<i_2<\ldots<i_m\}$. We associate with $J$ and $p$ the sequence $\Sigma(J,p)=(\sigma_1,\ldots,\sigma_{m+1})$, where $\sigma_j:=(\mathrm{sgn}(\gamma_{i_j}), \mathrm{sgn}(\gamma_{i_j}'))$ for $j=1,\ldots, m$, and $\sigma_{m+1}:= \mathrm{sgn}(\langle\gamma_{\infty},\alpha_p^\vee\rangle)$. We now present some properties of the sequence $\Sigma(J,p)$, which will be used later, and which reflect the combinatorics of admissible subsets, as discussed in \cite{LP1}.

\begin{proposition}\cite{LP1}\label{combprop} The sequence $\Sigma(J,p)$ has the following properties:
\begin{enumerate}
\item[(S1)] $\sigma_j\in\{(1,1),\,(-1,-1),\,(1,-1)\}$ for $j=1,\ldots, m$;
\item[(S2)] $j=0$ or $\sigma_j=(1,1)$ implies $\sigma_{j+1}\in\{(1,1),\,(1,-1),\,1\}$.
\end{enumerate}
\end{proposition}

The sequence $\Sigma(J,p)$ determines a continuous piecewise-linear function $g_{J,p}\::\:[0,m+\frac{1}{2}]\rightarrow{\mathbb R}$ as shown below. By a step $(h,k)$ of a function $f$ at $x=a$, we understand that $f(a+h)=f(a)+k$, and that $f$ is linear between $a$ and $a+h$. We set $g_{J,p}(0)=-\frac{1}{2}$ and, by scanning $\Sigma(J,p)$ from left to right while ignoring brackets, we impose the following condition: the $i$th entry $\pm 1$ corresponds to a step $(\frac{1}{2},\pm\frac{1}{2})$ of $g_{J,p}$ at $x=\frac{i-1}{2}$, respectively. 

\begin{proposition}\cite{LP1}\label{gap}
The function $g_{J,p}$ encodes the sequence $L(J,p)$ as follows:
\[l_{i_j}=g_{J,p}\!\left(j-\frac{1}{2}\right),\;j=1,\ldots,m\,,\;\;\;\; \mbox{and}\;\;\;\; l_\infty^p=g_{J,p}\!\left(m+\frac{1}{2}\right)\,.\]
\end{proposition}

\begin{example}\label{exgjp} Assume that the entries of $\Gamma(J)$ indexed by the elements of $I(J,p)$ are $(\alpha_p,-\alpha_p)$, $(-\alpha_p,-\alpha_p)$, $(\alpha_p,\alpha_p)$, $(\alpha_p,\alpha_p)$, $(\alpha_p,-\alpha_p)$, $(-\alpha_p,-\alpha_p)$, $(\alpha_p,-\alpha_p)$, $(\alpha_p,\alpha_p)$, in this order; also assume that ${\rm sgn}(\langle\gamma_\infty,\alpha_p^\vee\rangle)=1$. The graph of $g_{J,p}$ is shown in Figure~\ref{fig:walls}; this graph is separated into segments corresponding to the entries of the sequence $\Sigma(J,p)$. 
\begin{figure}[ht]
\mbox{\epsfig{file=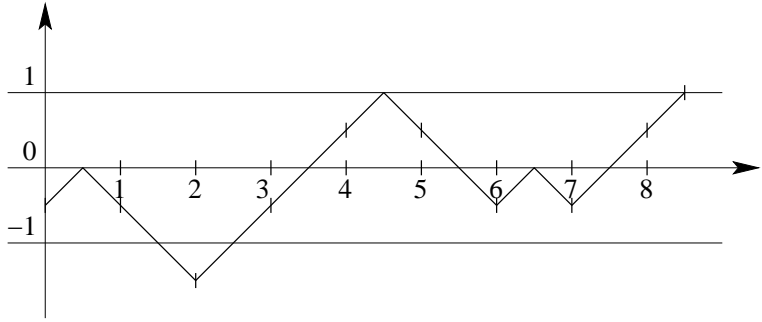}}
\caption{The graph of the function $g_{J,p}$ in Example \ref{exgjp}.}
\label{fig:walls}
\end{figure}
\end{example}

\subsection{Root operators}\label{r-ops} We now define partial operators known as {\em root operators} on the collection ${\mathcal A}(\Gamma)$ of admissible subsets corresponding to our fixed $\lambda$-chain. They are associated with a fixed simple root $\alpha_p$, and are traditionally denoted by $F_p$ (also called a lowering operator) and $E_p$ (also called a raising operator). The notation is the one introduced in the previous subsection.

We first consider  $F_p$ on the admissible subset $J$. 
 This is defined whenever $M(J,p)>0$. Let $m=m_F(J,p)$ be defined by
\[m_F(J,p):=\casetwo{\min\,\{i\in I(J,p)\mid l_i=M(J,p)\}}{\mbox{this set is nonempty}}{\infty}\]
Let $k=k_F(J,p)$ be the predecessor of $m$ in $I(J,p)\cup\{\infty\}$, which always exists. It turns out that $m\in J$ if $m\ne \infty$, but $k\not\in J$ (cf. Proposition \ref{propf} below). Finally, we set
\begin{equation}\label{deffp}
F_p(J):=(J\setminus\{m\})\cup \{k\}\,.
\end{equation}

\begin{proposition}\cite{LP1}\label{propf} Given the above setup, the following hold.
\begin{enumerate}
\item If $m\ne\infty$, then $\gamma_m'=-\gamma_m=-\alpha_p$. We also have $\gamma_k=\gamma_k'=\alpha_p$ and $l_k=M(J,p)-1$. 
\item We have ${\mu}(F_p(J))={\mu}(J)-\alpha_p\,.$
\item We have $w(F_p(J))=w(J)$ if $m\ne\infty$, and $w(F_p(J))=s_pw(J)$ otherwise.
\end{enumerate}
\end{proposition}

Let us now define a partial inverse $E_p$ to $F_p$. The operator $E_p$ is defined on the admissible subset $J$ whenever $M(J,p)>\langle \mu(J),\alpha_p^\vee\rangle$.  Let $k=k_E(J,p)$ be defined by
\[k_E(J,p):=\max\,\{i\in I(J,p)\mid l_i=M(J,p)\}\,;\]
the above set turns out to be always nonempty. Let $m=m_E(J,p)$ be the successor of $k$ in $I(J,p)\cup\{\infty\}$. It turns out that $k\in J$ but $m\not\in J$ (cf. Proposition \ref{prope} below). Finally, we set
\begin{equation}\label{defep}
E_p(J):=(J\setminus\{k\})\cup (\{m\}\setminus\{\infty\})\,.
\end{equation}

\begin{proposition}\cite{LP1}\label{prope} Given the above setup, the following hold.
\begin{enumerate}
\item We have $\gamma_k'=-\gamma_k=-\alpha_p$. If $m\ne\infty$, then $\gamma_m=\gamma_m'=-\alpha_p$, and $l_m=M(J,p)-1$. 
\item We have ${\mu}(E_p(J))={\mu}(J)+\alpha_p\,.$
\item We have $w(E_p(J))=w(J)$ if $m\ne\infty$, and $w(E_p(J))=s_pw(J)$ otherwise.
\end{enumerate}
\end{proposition}

Similarly to Kashiwara's operators (see Subsection \ref{cryst-inv}), the root operators above define a directed colored graph structure and a poset structure on the set ${\mathcal A}(\Gamma)$ of admissible subsets corresponding to a fixed $\lambda$-chain $\Gamma$. According to \cite[Proposition 6.9]{LP1}), the admissible subset $J_{\max}=\emptyset$ is the maximum of the poset ${\mathcal A}(\Gamma)$. The following result related to the special $\lambda$-chain in Proposition \ref{constr-lambdachain}, which we denote by $\Gamma^*$, was proved in \cite{LP1}.

\begin{theorem}\cite{LP1}\label{isom-cryst} The directed colored graph on  the set ${\mathcal A}(\Gamma^*)$ defined by the root operators is isomorphic to the crystal graph of the irreducible representation $V_\lambda$ with highest weight $\lambda$. Under this isomorphism, the weight of an admissible subset gives the weight space in which the corresponding element of the canonical basis lies. 
\end{theorem}

\section{Yang-Baxter Moves}\label{yb-moves}

In this section, we define the analog of Sch\"utzenberger's sliding algorithm in our model, which we call a {\em Yang-Baxter move}, for reasons explained below. We start with some results on dihedral subgroups of Weyl groups.

\subsection{Dihedral reflection subgroups} 
\begin{figure}[ht]
\mbox{\epsfig{file=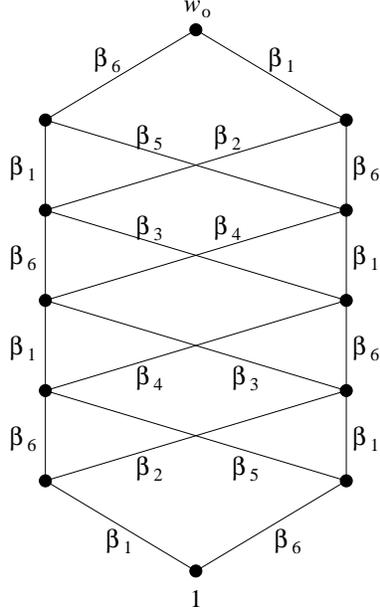}}
\caption{The Bruhat order on the Weyl group of type $G_2$.}
\label{bruhat}
\end{figure}

Let $\overline{W}$ be a dihedral Weyl group of order $2q$, that is, a Weyl group of type $A_1\times A_1$, $A_2$, $B_2$, or $G_2$ (with $q=2, \,3,\, 4,\, 6$, respectively). Let $\overline{\Phi}$ be the corresponding root system with simple roots $\alpha$, $\beta$. The sequence
\begin{equation}\label{reflorder}\beta_1:=\alpha,\;\;\beta_2:=s_\alpha(\beta),\;\;\beta_3:=s_\alpha s_\beta(\alpha),\;\;\ldots,\;\;\beta_{q-1}:=s_\beta(\alpha),\;\;\beta_q:=\beta\end{equation}
is a {\em reflection ordering} on the positive roots of $\overline{\Phi}$ (cf. \cite{Dyer}). The following Lemma describes the structure of $\overline{W}$ and its action on $\overline{\Phi}$. 
As an illustration, we present the Bruhat order on the Weyl group of type $G_2$ in Figure \ref{bruhat}. Here, as well as throughout this paper, we label a cover $v\lessdot vs_\gamma$ in Bruhat order by the corresponding root $\gamma$.

\begin{lemma}\label{perm-roots} {\rm (1)} If $i\le \frac{q+1}{2}$, then the reflection $s_{\beta_i}$ sends the roots $\beta_1,\dots,\beta_{i-1}$ to $-\beta_{2i-1},\dots,-\beta_{i+1}$, and the roots $\beta_{2i},\dots,\beta_q$ to $\beta_q,\dots,\beta_{2i}$, respectively. If $i> \frac{q+1}{2}$, then the reflection $s_{\beta_i}$ sends the roots $\beta_{i+1},\dots,\beta_{q}$ to $-\beta_{i-1},\dots,-\beta_{2i-q}$, and the roots $\beta_{1},\dots,\beta_{2i-q-1}$ to $\beta_{2i-q-1},\dots,\beta_{1}$, respectively.

{\rm (2)} Given $\overline{v}\in\overline{W}$ with $a:=\ell(\overline{v})<q$, consider its covers in Bruhat order by defining $\overline{\Phi}(\overline{v}):=\{j\in[q]\mid \ell(\overline{v}s_{\beta_j})=\ell(\overline{v})+1\}$. We have
\[\overline{\Phi}(\overline{v})=\casetwoex{\{1,q-a\}}{\overline{v}=\ldots s_\alpha s_\beta}{\{a+1,q\}}{\overline{v}=\ldots s_\beta s_\alpha}\]
\end{lemma}

With every pair of Weyl group elements $\overline{u}<\overline{w}$ in Bruhat order, we will associate a subset $J(\overline{u},\overline{w})$ of $[q]$ as follows. Let $a:=\ell(\overline{u})$ and $b:=\ell(\overline{w})$. Given $\delta\in\{\alpha,\beta\}$, we will use the notation 
\[\overline{W}_\delta:=\{\overline{v}\in \overline{W} \mid \ell(\overline{v}s_\delta)>\ell(\overline{v})\}\,,\;\;\;\;\overline{W}^\delta:=\overline{W}\setminus\overline{W}_\delta=\{\overline{v}\in \overline{W}\mid\ell(\overline{v}s_\delta)<\ell(\overline{v})\}\,.\]
\begin{enumerate}\item[]\begin{enumerate}
\item[{\bf Case 0:}] $\overline{u}=\overline{w}$. We let $J(\overline{u},\overline{u}):=\emptyset$.
\item[{\bf Case 1:}] $b-a=1$. We have the following disjoint subcases.
\begin{enumerate}
\item [{\bf Case 1.1:}] $\overline{u}\in \overline{W}_\alpha$, $\overline{w}\in\overline{W}^\alpha$, so $0\le a\le q-1$. We let $J(\overline{u},\overline{w}):=\{1\}$.
 \item [{\bf Case 1.2:}]  $\overline{u}\in \overline{W}^\beta$, $\overline{w}\in\overline{W}_\alpha$, so $0< a<q-1$. We let $J(\overline{u},\overline{w}):=\{q-a\}$.
\item [{\bf Case 1.3:}] $\overline{u}\in \overline{W}_\beta$, $\overline{w}\in\overline{W}^\beta$, so $0\le a\le q-1$. We let $J(\overline{u},\overline{w}):=\{q\}$.
\item [{\bf Case 1.4:}] $\overline{u}\in \overline{W}^\alpha$, $\overline{w}\in\overline{W}_\beta$, so $0< a<q-1$. We let $J(\overline{u},\overline{w}):=\{a+1\}$.
\end{enumerate}
\item[{\bf Case 2:}] $1<b-a<q$. We have the following disjoint subcases.
\begin{enumerate}
\item [{\bf Case 2.1:}] $\overline{u}\in \overline{W}_\alpha$, $\overline{w}\in\overline{W}_\beta$, so $0\le a<a+2\le b<q$. \\We let $J(\overline{u},\overline{w}):=\{1,a+2,a+3,\ldots,b\}$.
\item [{\bf Case 2.2:}] $\overline{u}\in \overline{W}^\beta$, $\overline{w}\in\overline{W}^\beta$, so $0< a<a+2\le b\le q$. \\We let $J(\overline{u},\overline{w}):=\{1,a+2,a+3,\ldots,b-1,q\}$.
\item [{\bf Case 2.3:}] $\overline{u}\in \overline{W}_\beta$, $\overline{w}\in\overline{W}_\alpha$, so $0\le a<a+2\le b<q$. \\We let $J(\overline{u},\overline{w}):=\{a+1,a+2,\ldots,b-1,q\}$.
\item [{\bf Case 2.4:}] $\overline{u}\in \overline{W}^\alpha$, $\overline{w}\in\overline{W}^\alpha$, so $0< a<a+2\le b\le q$. \\We let $J(\overline{u},\overline{w}):=\{a+1,a+2,\ldots,b\}$.
\end{enumerate}
\item[{\bf Case 3:}] $a=0$ and $b=q$, that is, $\overline{u}$ is the identity and $\overline{w}$ is the longest Weyl group element $\overline{w}_\circ$. In this case, we let $J:=[q]$.
\end{enumerate}
\end{enumerate}
In Case 2.2, if $b=a+2$ then the sequence $a+2,a+3,\ldots,b-1$ is considered empty.

Let $J(\overline{u},\overline{w}):=\{j_1<j_2<\ldots<j_{b-a}\}$. We use the notation $r_i:=s_{\beta_i}$, as above. In all cases above we have a unique saturated increasing chain in Bruhat order from $\overline{u}$ to $\overline{w}$ whose labels form a subsequence of (\ref{reflorder}); this chain is
\[\overline{u}\lessdot \overline{u}r_{{j_1}}\lessdot \overline{u}r_{{j_1}}r_{{j_2}}\lessdot\ldots\lessdot \overline{u}r_{{j_1}}\ldots r_{{j_{b-a}}}=\overline{w}\,.\]
Indeed, this can be easily checked based on Lemma \ref{perm-roots} (2). More generally, we have the  result below for an arbitrary Weyl group $W$ with a dihedral reflection subgroup $\overline{W}$ and corresponding root systems $\Phi\supseteq\overline{\Phi}$. The notation is the same as above. It is known that any element $w$ of $W$ can be written uniquely as $w=\floor{w}\overline{w}$, where $\floor{w}$ is the minimal representative of the left coset $w\overline{W}$, and $\overline{w}\in\overline{W}$.

\begin{proposition}\label{chain} For each pair of elements $u<w$ in the same (left) coset of $W$ modulo $\overline{W}$, we have a unique saturated increasing chain in Bruhat order from ${u}$ to ${w}$ whose labels form a subsequence of {\rm (\ref{reflorder})}; this chain is
\[{u}\lessdot {u}r_{{j_1}}\lessdot {u}r_{{j_1}}r_{{j_2}}\lessdot\ldots\lessdot {u}r_{{j_1}}\ldots r_{{j_{b-a}}}=w\,,\]
where $J(\overline{u},\overline{w}):=\{j_1<j_2<\ldots<j_{b-a}\}$.
\end{proposition}

This result can be easily deduced from the corresponding one for $W=\overline{W}$ via the following Lemma about cosets modulo dihedral reflection subgroups, which was discussed in \cite{bfpmbo}.

\begin{lemma}\cite{bfpmbo}
The Bruhat order on $\overline{W}$ (viewed as a Coxeter group with generators $s_\alpha$ and $s_\beta$) is isomorphic to the partial order on any coset $w\overline{W}$ (induced from the Bruhat order on $W$). The isomorphism is given by the map $\overline{w}\mapsto \floor{w}\overline{w}$. This statement can be rephrased by saying that, for any $\overline{w}\in \overline{W}$ and $\gamma\in\overline{\Phi}$, we have $\overline{w}<\overline{w}s_\gamma$ if and only if $\floor{w}\overline{w}<\floor{w}\overline{w}s_\gamma$. 
\end{lemma}

We obtain another reflection ordering by reversing the sequence (\ref{reflorder}). Let us denote the corresponding subset of $[q]$ by $J'(\overline{u},\overline{w})$. We are interested in passing from the chain between $u$ and $w$ compatible with the ordering (\ref{reflorder}) to the chain compatible with the reverse ordering. If we fix $a:=\ell(\overline{u})$ and $b:=\ell(\overline{w})$, we can realize the passage from $J(\overline{u},\overline{w})$ to $J'(\overline{u},\overline{w})$ via the involution $Y_{q,a,b}$ described below in each of the cases mentioned above.
\[\begin{array}{ll}
\!\!\!\!\!\!\!\!\!\!\!\mbox{\bf Case 0:} &\!\!\!\mbox{$\emptyset\leftrightarrow\emptyset$ if $a=b$}\,.\\[0.05in] 
\!\!\!\!\!\!\!\!\!\!\!\mbox{\bf Case 1.1:} &\!\!\!\mbox{$\{1\}\leftrightarrow\{q\}$ if $0\le a=b-1\le q-1$}\,.\\[0.05in]
\!\!\!\!\!\!\!\!\!\!\!\mbox{\bf Case 1.2:} &\!\!\!\mbox{$\{q-a\}\leftrightarrow\{a+1\}$ if $0<a=b-1<q-1$}\,.\\[0.05in]
\!\!\!\!\!\!\!\!\!\!\!\mbox{\bf Case 2.1:} &\!\!\!\mbox{$\{1,a+2,a+3,\ldots,b\}\leftrightarrow\{a+1,a+2,\ldots,b-1,q\}$ if $0\le a<a+2\le b<q$}\,.\\[0.05in]
\!\!\!\!\!\!\!\!\!\!\!\mbox{\bf Case 2.2:} &\!\!\!\mbox{$\{1,a+2,a+3,\ldots,b-1,q\}\leftrightarrow\{a+1,a+2,\ldots,b\}$ if $0< a<a+2\le b\le q$}\,.\\[0.05in]
\!\!\!\!\!\!\!\!\!\!\!\mbox{\bf Case 3:} &\!\!\!\mbox{$[q]\leftrightarrow[q]$ if $a=0$ and $b=q$}\,.
\end{array}\]

\subsection{Yang-Baxter moves and their properties} Let us now consider an index set 
\begin{equation}\label{indexset}I:=\{\overline{1}<\ldots<\overline{t}<1<\ldots<q<\overline{t+1}<\ldots<\overline{n}\}\,,\end{equation}
 and let $\overline{I}:=\{\overline{1},\ldots,\overline{n}\}$. Let $\Gamma=\{\beta_i\}_{i\in I}$ be a $\lambda$-chain, denote $r_i:=s_{\beta_i}$ as before, and let $\Gamma'=\{\beta_i'\}_{i\in I}$ be the sequence of roots defined by 
\begin{equation}\label{yblambdachain}\beta_i'=\casetwo{\beta_{q+1-i}}{i\in I\setminus\overline{I}}{\beta_i}\end{equation}
In other words, the sequence $\Gamma'$ is obtained from the $\lambda$-chain $\Gamma$ by reversing a certain segment. Now assume that $\{\beta_1,\ldots,\beta_q\}$ are the positive roots of a rank two root system $\overline{\Phi}$ (without repetition). Let $\overline{W}$ be the corresponding dihedral reflection subgroup of the Weyl group $W$. The following result is easily proved using the correspondence between $\lambda$-chains and reduced words for the affine Weyl group element $v_{-\lambda}$ mentioned in Definition  \ref{def:lambda-chain}; most importantly, we need to recall from the proof of \cite[Lemma 9.3]{LP} that the moves $\Gamma\rightarrow\Gamma'$ correspond to Coxeter moves (on the mentioned reduced words) in this context. 

\begin{proposition}\label{reversing} {\rm (1)} The sequence $\Gamma'$ is also a $\lambda$-chain, and the sequence $\beta_1,\ldots,\beta_q$ is a reflection ordering.

{\rm (2)} We can obtain any $\lambda$-chain for a fixed dominant weight $\lambda$ from any other $\lambda$-chain by moves of the form $\Gamma\rightarrow\Gamma'$. 
\end{proposition}

Let us now map the admissible subsets in ${\mathcal A}(\Gamma)$ to those in ${\mathcal A}(\Gamma')$. Given $J\in {\mathcal A}(\Gamma)$, let 
\begin{equation}\label{setup1}
\overline{J}:=J\cap\overline{I}\,,\;\;\;\;u:=w(J\cap\{\overline{1},\ldots,\overline{t}\})\,,\;\;\;\;\mbox{and}\;\;\;\;w:=w(J\cap(\{\overline{1},\ldots,\overline{t}\}\cup[q]))\,.\end{equation}
 Also let 
\begin{equation}\label{setup2}u=\floor{u}\overline{u}\,,\;\;\;\;w=\floor{w}\overline{w}\,,\;\;\;\;a:=\ell(\overline{u})\,,\;\;\;\;\mbox{and}\;\;\;\;b:=\ell(\overline{w})\,,\end{equation}
 as above. It is clear that we have a bijection $Y\::\:{\mathcal A}(\Gamma)\rightarrow {\mathcal A}(\Gamma')$ given by
\begin{equation}\label{setup3}Y(J):=\overline{J}\cup Y_{q,a,b}(J\setminus\overline{J})\,.\end{equation}
We call the moves $J\mapsto Y(J)$ {\em Yang-Baxter moves} (cf. the discussion following Theorem \ref{weightpres}). We say that they are of types 0, 1.1, 1.2, 2.1, 2.2, and 3 depending on the cases considered above in relation to the definition of the corresponding map $Y_{q,a,b}$; we also use the term type 1 (respectively 2) for types 1.1 or 1.2 (respectively 2.1 or 2.2). Clearly, a Yang-Baxter move preserves the Weyl group element $w(\,\cdot\,)$ associated to an admissible subset, that is, 
\begin{equation}\label{w-pres}w(Y(J))=w(J)\,.\end{equation}
In addition, Theorem \ref{weightpres} below holds.  

In order to prove the mentioned result, we need to recall some information from \cite{LP}. Consider the ring $K:=\Z[\Lambda/h]\otimes \Z[W]$, where $\Z[W]$ is the group algebra of the Weyl group $W$, and $\Z[\Lambda/h]$ is the group algebra of $\Lambda/h:=\{\lambda/h\mid \lambda\in\Lambda\}$ (i.e., of the weight lattice shrunk $h$ times, $h$ being the Coxeter number defined in Subsection \ref{rootsyst}). We define $\Z[\Lambda/h]$-linear operators $B_\alpha$ and $X^\lambda$ on $K$, where $\alpha$ is a positive root and $\lambda$ is a weight:
\[
B_\alpha : w\longmapsto\left\{\begin{array}{cl}
ws_\alpha & \textrm{if }  \ell(ws_\alpha) = \ell(w)+1 \\[.05in]
0 & \textrm{otherwise,}
\end{array}\right.\;\;\;\;\;
X^\lambda: w \mapsto e^{w(\lambda/h)} w.
\]
The following commutation relation will be needed:
\begin{equation}\label{commute}
B_\alpha \, X^\lambda = X^{s_\alpha(\lambda)} \, B_\alpha\,.
\end{equation}

\begin{theorem}\label{weightpres} The map $Y$ preserves the weight of an admissible subset. In other words, $\mu(Y(J))=\mu(J)$ for all admissible subsets $J$. 
\end{theorem}

\begin{proof}
Fix an admissible subset $J$ and, for each $i\in I$, let us set
\[Z_i:=\casetwo{B_{\beta_i}}{i\in J}{X^{\beta_i}}\]
We can calculate $\mu(J)$ as follows:
\begin{equation}\label{comp-weight}X^\rho Z_{\overline{n}}\ldots Z_{\overline{t+1}}Z_q\ldots Z_1 Z_{\overline{t}}\ldots Z_{\overline{1}}X^{-\rho}(1)=e^{\mu(J)} w(J)\,.\end{equation}
Indeed, let us denote the alcoves in the gallery $\gamma(J)$ by $A_i$ for $i\in\{\overline{0}\}\cup I$, and let us also consider the admissible folding $\Gamma(J)=(\{(\gamma_i,\gamma_i')\}_{i\in I},\gamma_\infty)$. Fix $i\in I$, and let $\zeta_i$ and $\zeta_{i'}$ be the central points of $A_i$ and $A_{i'}$, where $i'$ is the predecessor of $i$ in $\{\overline{0}\}\cup I$. Then, based on (\ref{gal-fol1}), we have
\[Z_i(e^\mu w)=\casetwo{e^\mu wr_i}{i\in J}{e^{\mu+\gamma_i/h} w=e^{\mu+\zeta_{i'}-\zeta_i}w}\]
Therefore, $Z_i\ldots Z_{\overline{1}}X^{-\rho}(1)=e^{-\zeta_i}w(J\cap\{j\in I\mid j\le i\})$. Finally, by (\ref{gal-fol1}), applying the last operator $X^\rho$ amounts to multiplying by $e^{\gamma_\infty/h}=e^{\zeta_{\overline{n}}+\mu(J)}$, where $\zeta_{\overline{n}}$ is the central point of $A_{\overline{n}}$. Denoting the operators $Z_i$ corresponding to $Y(J)$ by $Z_i'$, we will show that the compositions $Z_q\ldots Z_1$ and $Z_q'\ldots Z_1'$ coincide; hence, when plugging them into the left-hand side of (\ref{comp-weight}), we obtain the same result. 

The cases we now consider correspond to the types of the Yang-Baxter move $J\mapsto Y(J)$. If the set $J\cap[q]$ is empty or equal to $[q]$ (that is, we have a Yang-Baxter move of type 0 or 3), then we clearly have $Z_q\ldots Z_1=Z_q'\ldots Z_1'$.

{\bf Case 1:} $J\cap[q]=\{i\}$. We will show that the two compositions coincide, that is, 
\[X^{\beta_q}\ldots X^{\beta_{i+1}}B_{\beta_i}X^{\beta_{i-1}}\ldots X^{\beta_1}=X^{\beta_1}\ldots X^{\beta_{i-1}}B_{\beta_i}X^{\beta_{i+1}}\ldots X^{\beta_q}\,,\]
for $i=1,\ldots,q$. By commuting the two operators $B_{\beta_i}$ to the right, based on Lemma \ref{perm-roots} (1) and (\ref{commute}), both sides are equal to $X^{\beta_{2i}+\ldots+\beta_q}B_{\beta_i}$ if $i\le \frac{q+1}{2}$, and $X^{\beta_1+\ldots+\beta_{2i-q-1}}B_{\beta_i}$ if $i>\frac{q+1}{2}$. 

{\bf Case 2:} $J\cap[q]=\{1,a+2,a+3,\ldots,b\}$ or $J\cap[q]=\{a+1,a+2,\ldots,b-1,q\}$ with $0\le a<a+2\le b<q$. Assume that we have
\begin{align*}&Z_q\ldots Z_1=X^{\beta_q}\ldots X^{\beta_{b+1}}B_{\beta_b}\ldots B_{\beta_{a+2}}X^{\beta_{a+1}}\ldots X^{\beta_2}B_{\beta_1}\;\;\:\mbox{and}\\
&Z_q'\ldots Z_1'=B_{\beta_1}X^{\beta_2}\ldots X^{\beta_{q+1-b}}B_{\beta_{q+2-b}}\ldots B_{\beta_{q-a}}X^{\beta_{q+1-a}}\ldots X^{\beta_q}\,.
\end{align*}
 By commuting the two operators $B_{\beta_1}$ past the operators of the form $X^\mu$ to their left/right, based on Lemma \ref{perm-roots} (1) and (\ref{commute}), we obtain
\begin{align*}&Z_q\ldots Z_1=X^{\beta_q}\ldots X^{\beta_{b+1}}B_{\beta_b}\ldots B_{\beta_{a+2}}B_{\beta_1}X^{\beta_{q+1-a}}\ldots X^{\beta_q}\;\;\:\mbox{and}\\
&Z_q'\ldots Z_1'=X^{\beta_q}\ldots X^{\beta_{b+1}}B_{\beta_1}B_{\beta_{q+2-b}}\ldots B_{\beta_{q-a}}X^{\beta_{q+1-a}}\ldots X^{\beta_q}\,.
\end{align*}
 The case $J\cap[q]=\{1,a+2,a+3,\ldots,b-1,q\}$ or  $J\cap[q]=\{a+1,a+2,\ldots,b\}$ with $0<a<a+2\le b\le q$ is completely similar.
\end{proof}

We now explain the way in which the Yang-Baxter moves are related to the {\em Yang-Baxter equation}, which justifies the terminology. In \cite{LP}, we considered the operators $R_\alpha:=X^\rho(X^\alpha+B_\alpha)X^{-\rho}$ for $\alpha\in\Phi^+$; if $\alpha\in\Phi^-$, we defined $R_{\alpha}$ by setting $B_{\alpha}:=-B_{-\alpha}$. It was proved in \cite[Theorem 10.1]{LP} that the operators $\{R_\alpha \mid\alpha\in \Phi\}$ satisfy the Yang-Baxter equation in the sense of Cherednik \cite{Cher}. (In fact, the dual of $B_\alpha$ was used in \cite{LP}, but this does not affect the above result.) The main application of the operators $R_\alpha$ was to show that, given a $\lambda$-chain $\Gamma=(\beta_1,\ldots,\beta_n)$, we have
\begin{equation}\label{ybform}R_{\beta_n}\ldots R_{\beta_1}(1)=\sum_{J\in{\mathcal A}(\Gamma)}e^{\mu(J)}w(J)\,.\end{equation}
Due to the Yang-Baxter property, the right-hand side of the above formula does not change when we replace the $\lambda$-chain $\Gamma$ by $\Gamma'$, as defined above. The Yang-Baxter moves described above implement the passage from $\Gamma$ to $\Gamma'$ at the level of the individual terms in (\ref{ybform}).

Furthermore, let us note that Theorem \ref{weightpres} also follows by combining Proposition \ref{chain} with \cite[Theorem 10.1]{LP}, that was mentioned above. However, the proof of the latter theorem is based on an involved case by case check in \cite{bfpmbo}, while even the part of the proof in \cite{LP} is not transparent. By contrast, the proof of Theorem \ref{weightpres} presented here, based on making the map $Y$ explicit, is a direct and simple one.


\subsection{Yang-Baxter moves and root operators}\label{yb-root-ops} In this subsection, we present the main result related to Yang-Baxter moves.

 We start with a Lemma regarding the action of a root operator, which will be used several times below, and which is based on the combinatorics of admissible subsets discussed in Subsection \ref{comb-adm-seq}. As mentioned above, this combinatorics is best understood by graphing the piecewise-linear function $g_{J,p}$ associated to a simple root $\alpha_p$ and an admissible subset $J$. Let us also recall the definition of the set $I(J,p)$, of the sequence $L(J,p)$, and of the integer $M(J,p)$ in (\ref{ilmp}), as well as of the sequence $\Sigma(J,p)$. Finally, recall the definition of the positions $k_F(J,p)$ and $m_F(J,p)$ at the beginning of Subsection \ref{r-ops}, as well as Proposition \ref{propf}, which are all related to the root operator $F_p$.

\begin{lemma}\label{fp-act} Let $I(J,p)=\{i_1<i_2<\ldots<i_m\}$ and $\Gamma(J)=(\{(\gamma_i,\gamma_i')\}_{i\in I},\gamma_\infty)$. 

{\rm (1)} If we have
\[(\gamma_j,\gamma_j')=\casethreeexc{(-\alpha_p,-\alpha_p)}{j=i_c}{(\alpha_p,-\alpha_p)}{j=i_{c+1},i_{c+2},\ldots,i_{d-1}}{(\alpha_p,\alpha_p)}{j=i_d}\]
for some $1\le c<d\le m$, then $k_F(J,p)\ne i_d$.

{\rm (2)}  If we have
\[(\gamma_j,\gamma_j')=\casetwoexc{(\alpha_p,-\alpha_p)}{j=i_c,i_{c+1},\ldots,i_{d-1}}{(\alpha_p,\alpha_p)}{j=i_d}\]
for some $1\le c<d\le m$, then $m_F(J,p)\ne i_c$.
\end{lemma}

\begin{proof} Let $\Sigma(J,p)=(\sigma_1,\ldots,\sigma_{m+1})$. 

(1) Assume that $k_F(J,p)=i_d$. Then, by Proposition \ref{propf} (1), we have $M(J,p)=l_{i_d}+1$. We clearly have $l_{i_c}=l_{i_d}$. By Proposition \ref{combprop} (S2), we have $c>1$ and $\sigma_{c-1}\in\{(1,-1),\,(-1,-1)\}$. Therefore, by Proposition \ref{gap}, we have $l_{i_{c-1}}=l_{i_c}+1=M(J,p)$, which contradicts the definition of $k_F(J,p)$.  

(2) Assume that $m_F(J,p)=i_c$. Then $M(J,p)=l_{i_c}$. We clearly have $l_{i_c}=l_{i_d}$. By Proposition \ref{combprop} (S2), we have $\sigma_{d+1}\in\{(1,1),\,(1,-1),\,1\}$. Therefore, by Proposition \ref{gap}, we have $l_{i_{d+1}}=l_{i_d}+1=M(J,p)+1$ if $d<m$, or $l_\infty^p=l_{i_d}+1=M(J,p)+1$ if $d=m$. Both contradict the definition of $M(J,p)$.  
\end{proof}

\begin{theorem}\label{comm-f-y} The root operators commute with the Yang-Baxter moves, that is, a root operator $F_p$ is defined on an admissible subset $J$ if and only if it is defined on $Y(J)$ and we have
\[Y(F_p(J))=F_p(Y(J))\,.\]
\end{theorem}

\begin{proof} The setup is the one described above, particularly in (\ref{setup1})-(\ref{setup3}). Fix an admissible subset $J$ in ${\mathcal A}(\Gamma)$, and consider the corresponding admissible folding $\Gamma(J)=(\{(\gamma_i,\gamma_i')\}_{i\in I},\gamma_\infty)$. Let 
\[\overline{\Gamma}(J)=((\overline{\gamma}_1,\overline{\gamma}_1'),\ldots,(\overline{\gamma}_q,\overline{\gamma}_q')):=((\floor{u}^{-1}(\gamma_1),\floor{u}^{-1}(\gamma_1')),\ldots,(\floor{u}^{-1}(\gamma_q),\floor{u}^{-1}(\gamma_q')))\,.\]
Clearly, this sequence consists only of roots in $\overline{\Phi}$. We also consider restrictions of $\overline{\Gamma}(J)$ to subsets of consecutive elements $\{i,\,i+1,\,\ldots,\,j\}$ of the set $[q]$, which we denote by
\[\overline{\Gamma}(J)_i^j=((\overline{\gamma}_i,\overline{\gamma}_i'),(\overline{\gamma}_{i+1},\overline{\gamma}_{i+1}'),\ldots,(\overline{\gamma}_j,\overline{\gamma}_j'))\,.\]
Similar notation is used for any admissible subset, in particular for $Y(J)$.

Let $\alpha:=\beta_1$ and $\beta:=\beta_q$, as in (\ref{reflorder}). Note that the only indices $i\in[q]$ for which $\gamma_i$ or $-\gamma_i$ is a simple root are the ones for which $\overline{\gamma}_i$ belongs to $\{\pm\alpha,\,\pm\beta\}$. Indeed, if $\overline{\gamma}_i=c\alpha+d\beta$, then $\gamma_i=c\floor{u}(\alpha)+d\floor{u}(\beta)$, where $\floor{u}(\alpha)$ and $\floor{u}(\beta)$ are positive roots in $\Phi$ since $\ell(\floor{u}s_\alpha)>\ell(\floor{u})$ and $\ell(\floor{u}s_\beta)>\ell(\floor{u})$ (cf. \cite[Proposition~5.7]{Hum}). Hence, in order to compare the action of a root operator $F_p$ on $J$ and $Y(J)$, it is enough to consider the positions in $\overline{\Gamma}(J)$ and $\overline{\Gamma}(Y(J))$ in which the roots $\pm\alpha$ and $\pm\beta$ appear.

For simplicity, we denote the pairs of roots $(\gamma,\gamma)$ and $(\gamma,-\gamma)$ by $\gamma$ and $\pm\gamma$, respectively. It is also convenient to define
\[\delta(i)=\casetwoex{\alpha}{i \mbox{ odd}}{\beta}{i \mbox{ even}}\]

The cases we now consider, which depend on $\overline{u}$ and $\overline{w}$, are precisely the ones considered above in relation to the definition of the set $J(\overline{u},\overline{w})$; as discussed above, they give the type of the Yang-Baxter move $J\mapsto Y(J)$. The analysis below makes it clear that $F_p(J)$ and $F_p(Y(J))$ are both defined or undefined, so we assume that they are both defined whenever we mention them. If a root operator $F_p$ does not modify $J\cap[q]$ and $Y(J)\cap[q]$, then $F_p(J)$ and $F_p(Y(J))$ are clearly matched by a Yang-Baxter move of the same type as the one matching $J$ and $Y(J)$. Hence it suffices to assume that the root operator $F_p$ modifies $J\cap[q]$ or $Y(J)\cap[q]$. 

{\bf Case 0:} $J\cap[q]=\emptyset$. It is easy to see that $F_p(J)$ and $F_p(Y(J))$ are matched by a Yang-Baxter move of type 1.

{\bf Case 1.1:} $b-a=1$, $\overline{u}\in \overline{W}_\alpha$, $\overline{w}\in\overline{W}^\alpha$, so $0\le a\le q-1$ and $J\cap[q]=\{1\}$.  (Case 1.3 above is also treated here, since $Y(J)$ satisfies its conditions.) It is not hard to show that we have
\begin{align}\label{gj11}
&\overline{\Gamma}(J)_{a+1}^{a+2}=(-\delta(a+1),\delta(a+2))\,,\!\;\;\;\;\overline{\Gamma}(Y(J))_{a}^{a+1}=(-\delta(a+1),\delta(a+2))\,,\;\;\;\;\mbox{if $0<a<q-1\,,$}\\[0.02in]
\label{gj12}&\overline{\Gamma}(J)_{1}^{2}=(\pm\delta(1),\delta(2))\,,\;\;\;\;\;\;\;\;\;\;\;\;\;\;\;\;\;\;\;\;\overline{\Gamma}(Y(J))=(\delta(2),\ldots,\pm\delta(1))\,,\;\;\;\;\;\;\;\;\;\;\;\;\;\;\;\,\mbox{if $a=0\,,$}\\[0.02in]
\label{gj13}&\overline{\Gamma}(J)=(\pm\delta(q-1),\ldots,-\delta(q))\,,\,\;\;\;\;\;\:\overline{\Gamma}(Y(J))_{q-1}^{q}=(-\delta(q),\pm\delta(q-1))\,,\;\;\;\;\;\;\;\;\mbox{if $a=q-1\,.$}\end{align}

We present the proof of the first part of (\ref{gj11}), while the other facts can be proved similarly. Let $w_i$ be the element of $\overline{W}$ having length $i$ and the form $s_\alpha s_\beta \ldots$. We have 
\[\beta_i=w_{i-1}(\delta(i))\,,\;\;\;w_{i+1}=w_is_{\delta(i+1)}\,,\;\;\;\mbox{and}\;\;\;w(J\cap\{\overline{1},\ldots,\overline{t},1\})=us_\alpha=\floor{u}(\overline{u}s_{\alpha})=\floor{u}w_{a+1}^{-1}\,.\]
 Hence
\begin{equation}\label{gamma1}\begin{array}{l}\overline{\gamma}_{a+1}=w_{a+1}^{-1}(\beta_{a+1})=s_{\delta(a+1)}w_a^{-1}w_a(\delta(a+1))=-\delta(a+1)\,,\\[0.05in]
\overline{\gamma}_{a+2}=w_{a+1}^{-1}(\beta_{a+2})=w_{a+1}^{-1}w_{a+1}(\delta(a+2))=\delta(a+2)\,.\end{array}\end{equation} 

Note that the roots $\pm\alpha$ and $\pm\beta$ do not appear in other positions in $\overline{\Gamma}(J)$ and $\overline{\Gamma}(Y(J))$ beside the ones indicated in (\ref{gj11})-(\ref{gj13}). For instance, one can show this for the first part of (\ref{gj11}) by an argument completely similar to the one used in Case 2.1 below relative to the first part of (\ref{gj21}).

In (\ref{gj11}), the root operator $F_p$ must insert $a+2$ into $J$ and $a+1$ into $Y(J)$; hence $F_p(J)$ and $F_p(Y(J))$ are matched by a Yang-Baxter move of type 2.1 if $a<q-2$ (more precisely, $\{1,a+2\}\leftrightarrow\{a+1,q\}$) and by a move of type 2.2 if $a=q-2$ (more precisely, $\{1,q\}\leftrightarrow\{q-1,q\}$). In (\ref{gj12}), the root operator $F_p$ must either remove $1$ from $J$ and $q$ from $Y(J)$, or insert 2 into $J$ and 1 into $Y(J)$; hence $F_p(J)$ and $F_p(Y(J))$ are matched by a Yang-Baxter move of type 0, 2.1 (more precisely, $\{1,2\}\leftrightarrow\{1,q\}$), or 3 (this case is the analog of the previous one for $\overline{\Phi}$ of type $A_1\times A_1$). In (\ref{gj13}), the root operator $F_p$ must remove $1$ from $J$ and $q$ from $Y(J)$; hence $F_p(J)$ and $F_p(Y(J))$ are matched by a Yang-Baxter move of type 0.

{\bf Case 1.2:} $b-a=1$, $\overline{u}\in \overline{W}^\beta$, $\overline{w}\in\overline{W}_\alpha$, so $0< a<q-1$ and $J\cap[q]=\{q-a\}$. (Case 1.4 above is also treated here, since $Y(J)$ satisfies its conditions.) In a similar way to (\ref{gj11})-(\ref{gj13}), we can prove that we have
\begin{equation}\label{gj14}
\overline{\Gamma}(J)_{q-a}^{q-a}=(\pm\delta(a))\,,\;\;\;\;\overline{\Gamma}(Y(J))_{a}^{a+2}=(-\delta(a-1),\pm\delta(a),\delta(a+1))\,.
\end{equation}

As in the previous case, one can easily show that the roots $\pm\alpha$ and $\pm\beta$ do not appear in other positions in $\overline{\Gamma}(J)$ and $\overline{\Gamma}(Y(J))$ beside the ones indicated in (\ref{gj14}).

In (\ref{gj14}), the root operator $F_p$ must remove $q-a$ from $J$ and $a+1$ from $Y(J)$. Hence $F_p(J)$ and $F_p(Y(J))$ are matched by a Yang-Baxter move of type 0. Note that $F_p$ cannot insert $a+2$ into $Y(J)$, by Lemma \ref{fp-act} (1). 

{\bf Case 2.1:} $1<b-a<q$, $\overline{u}\in \overline{W}_\alpha$, $\overline{w}\in\overline{W}_\beta$, so $0\le a<a+2\le b<q$ and $J\cap[q]=\{1,a+2,a+3,\ldots,b\}$. (Case 2.3 above is also treated here, since $Y(J)$ satisfies its conditions.) We start by showing that we have
\begin{equation}\begin{array}{l}\label{gj21}\overline{\Gamma}(J)_{a+1}^{b+1}=(-\delta(a+1),\pm\delta(a+2),\pm\delta(a+3),\ldots,\pm\delta(b),\delta(b+1))\\[0.05in]
\overline{\Gamma}(Y(J))_{a}^{b}=(-\delta(a+1),\pm\delta(a+2),\pm\delta(a+3),\ldots,\pm\delta(b),\delta(b+1))\end{array}\;\;\;\;\mbox{if $a>0\,,$}\end{equation}
as well as
\begin{equation}\begin{array}{l}\label{gj22}\overline{\Gamma}(J)_{1}^{b+1}=(\pm\delta(1),\pm\delta(2),\ldots,\pm\delta(b),\delta(b+1))\\[0.05in]
\overline{\Gamma}(Y(J))_{1}^{b}=(\pm\delta(2),\pm\delta(3),\ldots,\pm\delta(b),\delta(b+1))\end{array}\;\;\;\;\mbox{if $a=0\,.$}\end{equation}

We present the proof of the first part of (\ref{gj21}), while the other facts can be proved similarly. The roots $\overline{\gamma}_{a+1}$ and $\overline{\gamma}_{a+2}$ can be computed as in (\ref{gamma1}). For $i=a+3,\ldots,b+1$, we calculate based on Lemma \ref{perm-roots} (1) and (\ref{gamma1}):
\begin{equation}\label{gamma2}\begin{array}{ll}\overline{\gamma}_i&\!\!\!\!=w_{a+1}^{-1}s_{\beta_{a+2}}\ldots s_{\beta_{i-1}}(\beta_{i})=w_{a+1}^{-1}s_{\beta_{a+2}}\ldots s_{\beta_{i-2}}(-\beta_{i-2})=w_{a+1}^{-1}s_{\beta_{a+2}}\ldots s_{\beta_{i-3}}(\beta_{i-2})\\[0.05in]&\!\!\!\!=\ldots=\casetwoex{w_{a+1}^{-1}(-\beta_{a+1})=\delta(a+1)=\delta(i)}{i-a \mbox{ odd}}{w_{a+1}^{-1}(\beta_{a+2})=\delta(a+2)=\delta(i)}{i-a \mbox{ even}}\end{array}\end{equation}

Let us also note that the roots $\pm\alpha$ and $\pm\beta$ do not appear in other positions in $\overline{\Gamma}(J)$ and $\overline{\Gamma}(Y(J))$ beside the ones indicated in (\ref{gj21})-(\ref{gj22}). For instance, in the first part of (\ref{gj21}), we have $\overline{\gamma}_i=\pm w_{a+1}^{-1}(\beta_{i})\not\in\{\pm\alpha,\,\pm\beta\}$ for $i=1,\ldots,a$, due to (\ref{gamma1}). Similarly, in the same case, we have $\overline{\gamma}_j=w_{a+1}^{-1}s_{\beta_{a+2}}\ldots s_{\beta_{b}}(\beta_{j})\not\in\{\pm\alpha,\,\pm\beta\}$ for $j=b+2,\ldots,q$, based on (\ref{gamma2}) for $i=b,b+1$.

One way in which the operator $F_p$ can act on $J\cap[q]$ and $Y(J)\cap[q]$ is to insert $b+1$ into $J$ and $b$ into $Y(J)$. This can happen both in (\ref{gj21}) and in (\ref{gj22}), but in the former case only if $b-a$ is odd, by Lemma \ref{fp-act} (1). Hence $F_p(J)$ and $F_p(Y(J))$ are matched by a move of the form $\{1,a+2,a+3,\ldots,b,b+1\}\leftrightarrow\{a+1,a+2,\ldots,b-1,b,q\}$. This is a Yang-Baxter move of type 2.1 if $b<q-1$, of type 2.2 if $b=q-1,\;a>0$, and of type 3 if $b=q-1,\;a=0$. 

Finally, we consider the case when $F_p$ removes certain elements from $J\cap[q]$ and $Y(J)\cap[q]$. Let us first concentrate on the case $a>0$. Then $F_p$ must remove $a+2$ from $J$ and $a+1$ from $Y(J)$, but this can only happen if $b-a$ is even, by Lemma \ref{fp-act} (2). Thus $F_p(J)$ and $F_p(Y(J))$ are matched by a move of the form $\{1,a+3,\ldots,b\}\leftrightarrow\{a+2,\ldots,b-1,q\}$. This is a Yang-Baxter move of type 2.1 if $a+3\le b$, and of type 1.1 if $b=a+2$. Now let us turn to the case $a=0$. If $b$ is odd, then $F_p$ must remove 1 from $J$ and 2 from $Y(J)$, by Lemma \ref{fp-act} (2). Thus $F_p(J)$ and $F_p(Y(J))$ are matched by a move of the form $\{2,3,\ldots,b\}\leftrightarrow\{1,3,4,\ldots,b-1,q\}$, which is a Yang-Baxter move of type 2.2. If $b$ is even, then $F_p$ must remove 2 from $J$ and 1 from $Y(J)$, by Lemma \ref{fp-act} (2).  Thus $F_p(J)$ and $F_p(Y(J))$ are matched by a move of the form $\{1,3,4,\ldots,b\}\leftrightarrow\{2,3,\ldots,b-1,q\}$, which is a Yang-Baxter move of type 2.1 if $b>2$ and of type 1.1 if $b=2$. 

{\bf Case 2.2:}  $1<b-a<q$, $\overline{u}\in \overline{W}^\beta$, $\overline{w}\in\overline{W}^\beta$, so $0< a<a+2\le b\le q$ and $J\cap[q]=\{1,a+2,a+3,\ldots,b-1,q\}$. (Case 2.4 above is also treated here, since $Y(J)$ satisfies its conditions.) In a similar way to (\ref{gj21})-(\ref{gj22}), we can prove that we have
\begin{equation}\begin{array}{l}\label{gj23}\overline{\Gamma}(J)_{a+1}^{b}=(-\delta(a+1),\pm\delta(a+2),\pm\delta(a+3),\ldots,\pm\delta(b-1),\delta(b))\\[0.05in]
\overline{\Gamma}(Y(J))_{a}^{b+1}=(-\delta(a+1),\pm\delta(a+2),\pm\delta(a+3),\ldots,\pm\delta(b),\pm\delta(b+1),\delta(b+2))\end{array}\;\;\;\;\mbox{if $b<q\,,$}\end{equation}
as well as
\begin{equation}\label{gj24}\begin{array}{l}\overline{\Gamma}(J)_{a+1}^{q}=(-\delta(a+1),\pm\delta(a+2),\pm\delta(a+3),\ldots,\pm\delta(q-1),\pm\delta(q))\\[0.05in]
\overline{\Gamma}(Y(J))_{a}^{q}=(-\delta(a+1),\pm\delta(a+2),\pm\delta(a+3),\ldots,\pm\delta(q),\pm\delta(q+1))\end{array}\;\;\;\;\mbox{if $b=q\,.$}\end{equation}

As in the previous cases, one can easily show that the roots $\pm\alpha$ and $\pm\beta$ do not appear in other positions in $\overline{\Gamma}(J)$ and $\overline{\Gamma}(Y(J))$ beside the ones indicated in (\ref{gj23})-(\ref{gj24}).

One way in which the operator $F_p$ can act on $J\cap[q]$ and $Y(J)\cap[q]$ is to insert $b$ into $J$ and $b+1$ into $Y(J)$. This can happen in (\ref{gj23}), but only if $b-a$ is even, by Lemma \ref{fp-act} (1). Hence $F_p(J)$ and $F_p(Y(J))$ are matched by a move of the form $\{1,a+2,a+3,\ldots,b-1,b,q\}\leftrightarrow\{a+1,a+2,\ldots,b,b+1\}$. This is always a Yang-Baxter move of type 2.2. 

Finally, we consider the case when $F_p$ removes certain elements from $J\cap[q]$ and $Y(J)\cap[q]$. Then it must remove $a+2$ from $J$ and $a+1$ from $Y(J)$. This can happen both in (\ref{gj23}) and in (\ref{gj24}), but in the former case only if $b-a$ is odd, by Lemma \ref{fp-act} (2). Hence $F_p(J)$ and $F_p(Y(J))$ are matched by a move of the form $\{1,a+3,a+4,\ldots,b-1,q\}\leftrightarrow\{a+2,a+3,\ldots,b\}$. This is a Yang-Baxter move of type 2.2 with the exception of the case $b=q,\,a=q-2$, when it is of type 1.1. 

{\bf Case 3:} $a=0$ and $b=q$, so $[q]\subseteq J$. In this case we have
\[\overline{\Gamma}(J)=(\pm\delta(1),\pm\delta(2),\ldots)\,,\;\;\;\;  \overline{\Gamma}(Y(J))=(\pm\delta(2),\pm\delta(3),\ldots)\,.\]
Our root operator $F_p$ must either remove 1 from $J$ and 2 from $Y(J)$, or 2 from $J$ and 1 from $Y(J)$. Hence $F_p(J)$ and $F_p(Y(J))$ are matched by a Yang-Baxter move of type 2.2 (more precisely, $\{1,3,4,\ldots,q-1,q\}\leftrightarrow\{2,3,\ldots,q\}$) or 1.1 (this case is the analog of the previous one for $\overline{\Phi}$ of type $A_1\times A_1$). \end{proof}

Theorem \ref{comm-f-y} asserts that the map $Y$ above is an isomorphism between ${\mathcal A}(\Gamma)$ and ${\mathcal A}(\Gamma')$ as directed colored graphs. Given two arbitrary $\lambda$-chains $\Gamma$ and $\Gamma'$, we know from Proposition \ref{reversing} (2) that they can be related by a sequence of $\lambda$-chains $\Gamma=\Gamma_0,\,\Gamma_1,\,\ldots,\, \Gamma_m=\Gamma'$ to which correspond Yang-Baxter moves $Y_1,\,\ldots,\,Y_m$. Hence the composition $Y_m\ldots Y_1$ is an isomorphism between ${\mathcal A}(\Gamma)$ and ${\mathcal A}(\Gamma')$ as directed colored graphs. Since every directed graph ${\mathcal A}(\Gamma)$ has a unique source (cf. \cite[Proposition 6.9]{LP1}), its automorphism group as a directed colored graph consists only of the identity. Thus, we have the following corollary of Theorem \ref{comm-f-y}.

\begin{corollary}\label{isom-graphs} Given two arbitrary $\lambda$-chains $\Gamma$ and $\Gamma'$, the directed colored graph structures on ${\mathcal A}(\Gamma)$ and ${\mathcal A}(\Gamma')$ are isomorphic. This isomorphism is unique and, therefore, is given by the composition of Yang-Baxter moves corresponding to {\em any} sequence of $\lambda$-chains relating $\Gamma$ and $\Gamma'$. 
\end{corollary}

We have given a transparent combinatorial explanation for the independence of the directed colored graph defined by our root operators from the chosen $\lambda$-chain. Similarly, it was proved in \cite{Li2} that the directed colored graph structure on Littelmann paths generated by the corresponding root operators is independent of the initial path. However, this proof, which is based on continuous arguments, is less transparent. 

Based on Corollary \ref{isom-graphs}, Theorem \ref{isom-cryst} immediately leads to its generalization below. 

\begin{corollary}\label{all-cryst}
Given {\em any} $\lambda$-chain $\Gamma$, the directed colored graph on  the set ${\mathcal A}(\Gamma)$ defined by the root operators is isomorphic to the crystal graph of the irreducible representation $V_\lambda$ with highest weight $\lambda$. Under this isomorphism, the weight of an admissible subset gives the weight space in which the corresponding element of the canonical basis lies.
\end{corollary}

Based on Corollary \ref{all-cryst}, we will now identify the elements of the canonical basis with the corresponding admissible subsets.

\begin{remark}
We suggest that  root operators and Yang-Baxter moves would be able to explain the whole combinatorics of our model. Note the analogy with type $A$, where we have left strings and right strings, defined via root operators and jeu de taquin, respectively (cf. \cite{lasdcg}). 
\end{remark}

Define an action of a simple reflection $s_p$ on an admissible subset $J$ by 
\begin{equation}\label{w-action}s_p(J):=F_p^{\langle\mu(J),\alpha_p^\vee\rangle}(J)\,.\end{equation}
Up to the isomorphism in Corollary \ref{all-cryst}, this action coincides with the one on crystals defined by Kashiwara in \cite{kascbm} and \cite[Theorem 11.1]{kascb}; hence it leads to an action of the full Weyl group $W$. 

\begin{corollary} Equation {\rm (\ref{w-action})} defines a $W$-action on admissible subsets. We have $\mu(w(J))=w(\mu(J))$ for all $w$ in $W$ and all admissible subsets $J$.
\end{corollary}

\section{Lusztig's Involution}\label{sch-inv}

In this section, we present an explicit description of the involution $\eta_\lambda$ in Subsection \ref{cryst-inv} in the spirit of Sch\"utzenberger's evacuation. We will show that the role of jeu de taquin in the definition of the evacuation map is played by the Yang-Baxter moves.

\subsection{Reversing $\lambda$-chains and admissible subsets}\label{rev-lam-ch-j} Throughout the remainder of this paper, we fix an index set $I:=\{\overline{1}<\ldots<\overline{q}<1<\ldots<n\}$ and a $\lambda$-chain $\Gamma=\{\beta_i\}_{i\in I}$ such that $l_i^\emptyset=0$ if and only if $i\in \overline{I}:=\{\overline{1}<\ldots<\overline{q}\}$. In other words, the second occurence of a root can never be before the first occurence of another root. We will also write $\Gamma:=(\beta_{\overline{1}},\ldots,\beta_{\overline{q}},\beta_1,\ldots,\beta_n)$. Let us recall the notation $r_i:=s_{\beta_i}$ for $i\in I$. 

Given a Weyl group element $w$, we denote by $\floor{w}$ and $\ceil{w}$ the minimal and the maximal representatives of the coset $wW_\lambda$, respectively (where $W_\lambda$ is the stabilizer of the weight $\lambda$). Let $w_\circ^\lambda$ be the longest element of $W_\lambda$. Based on the discussion in Subsection \ref{lamch}, it is easy to see that we have the saturated increasing chain in Bruhat order 
\[1\lessdot r_{\overline{1}}\lessdot r_{\overline{1}}r_{\overline{2}}\lessdot\ldots\lessdot r_{\overline{1}}\ldots r_{\overline{q}}\]
 from 1 to $\floor{w_\circ}=w_\circ w_\circ^\lambda$. Hence the set $J_{\min}:=\overline{I}$ is an admissible subset. 

\begin{proposition}
The admissible subset $J_{\min}$ is the minimum of the poset ${\mathcal A}(\Gamma)$. 
\end{proposition} 

\begin{proof}
It suffices to show that, for any admissible subset $J\ne J_{\min}$, there exists $p\in[r]$ such that $M(J,p)>0$; in other words, the root operator $F_p$ is defined on $J$. Indeed, given such $J$, let ${j}=\min\,\overline{I}\setminus J$, which exists. Let $\Gamma(J)=(\{(\gamma_i,\gamma_i')\}_{i\in I},\gamma_\infty)$. It follows from definitions that $\gamma_j=\gamma_j'$ is a simple root $\alpha_p$. Proposition \ref{combprop} (S2) then implies $M(J,p)>0$.
\end{proof}

\begin{definition}\label{defkeys} Let $J$ be an admissible subset. Let $J\cap \overline{I}=\{\overline{j}_1<\ldots<\overline{j}_a\}$ and $J\setminus \overline{I}=\{j_1<\ldots<j_s\}$. The {\em initial key} $\kappa_0(J)$ and the {\em final key} $\kappa_1(J)$ of $J$ are the Weyl group elements defined by
\[\kappa_0(J):=r_{\overline{j}_1}\ldots r_{\overline{j}_a}\,,\;\;\;\;\;\kappa_1(J):=w(J)=\kappa_0(J)r_{j_1}\ldots r_{j_s}\,.\]
\end{definition}

\begin{remark}\label{remkeys} The keys $\kappa_0(J)$ and $\kappa_1(J)$ are the generalizations of the left and right keys of a semistandard Young tableau \cite{lasksb}, respectively. They are interchanged by Lusztig's involution (cf. Corollary \ref{conjkeys}) and are related to the Demazure character formula in Theorem \ref{demform}. Now recall the bijection in \cite[Section 9]{LP1} between LS chains (in the orbit of $-\lambda$) and admissible subsets for the special $\lambda$-chain. It is not hard to show that $\kappa_0(J)(-\lambda)$ and $\kappa_1(J)(-\lambda)$ are the initial and the final directions of the LS chain associated to $J$, respectively. If, instead, we use LS chains in the orbit of $\lambda$ (as we usually do), then $\kappa_0(J)(\lambda)$ and $\kappa_1(J)(\lambda)$ are the final and the initial directions of the corresponding LS chain, respectively.  
\end{remark}

We associate with our fixed $\lambda$-chain $\Gamma$ another sequence $\Gamma^{\rm rev}:=\{\beta_i'\}_{i\in I}$ by
\[\beta_i':=\casetwo{\beta_i}{i\in\overline{I}}{w_\circ^\lambda(\beta_{n+1-i})}\]
In other words, we have
\begin{equation}\label{defgrev}\Gamma^{\rm rev}=(\beta_{\overline{1}},\ldots,\beta_{\overline{q}},w_\circ^\lambda(\beta_n),w_\circ^\lambda(\beta_{n-1}),\ldots,w_\circ^\lambda(\beta_1))\,.\end{equation}

\begin{proposition}\label{gr-lamch} $\Gamma^{\rm rev}$ is a $\lambda$-chain. 
\end{proposition}

\begin{proof}
Note first that $w_\circ^\lambda$ permutes the roots in $\Phi^+\setminus\Phi_\lambda$, because so does any simple reflection in $W_\lambda$; here $\Phi_\lambda$ is the parabolic subroot system corresponding to $W_\lambda$. Therefore, since the $\lambda$-chain $\Gamma$ consists only of roots in $\Phi^+\setminus\Phi_\lambda$, so does $\Gamma^{\rm rev}$. 

We use the characterization of $\lambda$-chains in Theorem \ref{equivdef} (c). We observe first that the number of occurences of any positive root $\alpha$ in $\Gamma^{\rm rev}$ is $\langle\lambda,\alpha^\vee\rangle$. Indeed, if $\alpha\in \Phi^+\setminus\Phi_\lambda$, we have $\langle\lambda,w_\circ^\lambda(\alpha)^\vee\rangle=\langle\lambda,\alpha^\vee\rangle$. 

Let us now fix three positive roots $\alpha,\beta,\gamma$ such that $\gamma^\vee=\alpha^\vee+\beta^\vee$. Assume first that $\langle\lambda,\alpha^\vee\rangle$ and $\langle\lambda,\beta^\vee\rangle$ are both nonzero. Consider the subsequence of $\{\beta_i\}_{i\in I\setminus \overline{I}}$ consisting of $w_\circ^\lambda(\alpha)$, $w_\circ^\lambda(\beta)$, and $w_\circ^\lambda(\gamma)$. This starts with  $w_\circ^\lambda(\gamma)$, and continues with a concatenation of pairs $(w_\circ^\lambda(\alpha), w_\circ^\lambda(\gamma))$ and $(w_\circ^\lambda(\beta),w_\circ^\lambda(\gamma))$. Hence, the subsequence of $\{\beta_i'\}_{i\in I\setminus \overline{I}}$ consisting of $\alpha$, $\beta$, and $\gamma$ starts with $\gamma$ and continues with a concatenation of pairs $(\alpha,\gamma)$ and $(\beta,\gamma)$. Also, the subsequence of $\{\beta_i'\}_{i\in \overline{I}}$ consisting of $\alpha$, $\beta$, and $\gamma$ is either $(\alpha,\gamma,\beta)$ or $(\beta,\gamma,\alpha)$. 

Now assume that $\langle\lambda,\alpha^\vee\rangle=0$ and $\langle\lambda,\beta^\vee\rangle>0$. 
The subsequence of $\{\beta_i\}_{i\in I\setminus \overline{I}}$ consisting of $w_\circ^\lambda(\alpha)$, $w_\circ^\lambda(\beta)$, and $w_\circ^\lambda(\gamma)$ is a concatenation of pairs $(w_\circ^\lambda(\beta),w_\circ^\lambda(\gamma))$. Hence, the subsequence of $\{\beta_i'\}_{i\in I\setminus \overline{I}}$ consisting of $\alpha$, $\beta$, and $\gamma$ is a concatenation of pairs $(\beta,\gamma)$. Also, the subsequence of $\{\beta_i'\}_{i\in \overline{I}}$ consisting of $\alpha$, $\beta$, and $\gamma$ is  $(\beta,\gamma)$. 
\end{proof}

Let $r_i':=s_{\beta_i'}$ for $i\in I$. Fix an admissible subset 
\begin{equation}\label{defj}
J=\{\overline{j}_1<\ldots<\overline{j}_a<j_1<\ldots<j_s\} \end{equation}
in ${\mathcal A}(\Gamma)$, where $\{\overline{j}_1<\ldots<\overline{j}_a\}\subseteq\overline{I}$ and $\{j_1<\ldots<j_s\}\subseteq I\setminus\overline{I}$. Let $u:=\kappa_0(J)$ and $w:=\kappa_1(J)$. We have the increasing saturated chain
\begin{equation}\label{satch}
1\lessdot r_{\overline{j}_1}\lessdot r_{\overline{j}_1}r_{\overline{j}_2}\lessdot\ldots\lessdot r_{\overline{j}_1}\ldots r_{\overline{j}_a}=u\lessdot ur_{j_1}\lessdot ur_{j_1}r_{j_2}\lessdot\ldots\lessdot ur_{j_1}\ldots r_{j_s}=w\,.
\end{equation}
According to \cite{Dyer}, there is a unique saturated increasing chain in Bruhat order of the form
\[1\lessdot r_{\overline{k}_1}'\lessdot r_{\overline{k}_1}'r_{\overline{k}_2}'\lessdot\ldots\lessdot r_{\overline{k}_1}'\ldots r_{\overline{k}_b}'=\lfloor w_\circ w\rfloor=w_\circ w w_\circ^\lambda\,,\]
where $\{\overline{k}_1<\overline{k}_2<\ldots<\overline{k}_b\}\subseteq\overline{I}$. Define
\begin{equation}\label{defjrev}J^{\rm rev}:=\{\overline{k}_1<\ldots<\overline{k}_b<k_1<\ldots<k_s\}\,,\end{equation}
where $k_i:=n+1-j_{s+1-i}$ for $i=1,\ldots,s$. Note that $\beta_{k_i}'=w_\circ^\lambda(\beta_{j_{s+1-i}})$ for $i=1,\ldots,s$.

\begin{proposition}\label{jrev} $J^{\rm rev}$ is an admissible subset in ${\mathcal A}(\Gamma^{\rm rev})$. We have
\begin{equation}\label{keys-rev}\kappa_0(J^{\rm rev})=\floor{w_\circ \kappa_1(J)}\,,\;\;\;\;\kappa_1(J^{\rm rev})=\floor{w_\circ \kappa_0(J)}\,,\end{equation}
as well as $(J^{\rm rev})^{\rm rev}=J$.
\end{proposition}

\begin{proof}
We have $r_{k_i}'=w_\circ^\lambda r_{j_{s+1-i}}w_\circ^\lambda$. Therefore, according to (\ref{satch}), we have the saturated increasing chain
\begin{align*}&\lfloor w_\circ w\rfloor=w_\circ w w_\circ^\lambda \lessdot w_\circ w w_\circ^\lambda r_{k_1}'=w_\circ w r_{j_s}w_\circ^\lambda\lessdot w_\circ w w_\circ^\lambda r_{k_1}'r_{k_2}'=w_\circ w r_{j_s}r_{j_{s-1}}w_\circ^\lambda\lessdot\\ &\ldots\lessdot w_\circ w w_\circ^\lambda r_{k_1}'\ldots r_{k_s}'=w_\circ w r_{j_s} \ldots r_{j_1}w_\circ^\lambda=w_\circ u w_\circ^\lambda=\lfloor w_\circ u\rfloor\,.\end{align*}
This completes the proof of (\ref{keys-rev}), which then easily implies the last statement.
\end{proof}

 We now present a direct way to obtain the gallery $\gamma(J^{\rm rev})$ from $\gamma(J)$. Let us write 
\[\gamma(J)=(F_{\overline{0}},A_{\overline{0}},F_{\overline{1}},\ldots,F_{\overline{q}},A_{0},F_1,A_1,\ldots,A_n,F_{\infty})\,;\]
the corresponding augmented index set is $\{\overline{0}<\overline{1}<\ldots<\overline{q}=0<1<\ldots<n<\infty\}$. Let $\mu:=-\mu(J)$, that is, $F_{\infty}=\{\mu\}$. Now define another gallery in the following way:
\[\gamma^\omega:=(F_{\overline{0}}',A_{\overline{0}}',F_{\overline{1}}',\ldots,F_{\overline{q}}',A_{0}',F_1',A_1',\ldots,A_n',F_{\infty}')\,.\]
The notation is as follows:
\begin{itemize}
\item $\omega$ is the map on $\hR$ defined by $x\mapsto -w_\circ(x-\mu)$;
\item $A_i':=\omega(A_{n-i})$ for $i=0,\ldots,n$, $F_i':=\omega(F_{n+1-i})$ for $i=1,\ldots,n$, and $F_{\infty}'=\{w_\circ(\mu)\}$;
\item $(F_{\overline{0}}',A_{\overline{0}}',F_{\overline{1}}',\ldots,F_{\overline{q}}')$ is the initial segment of the gallery $\gamma(J^{\rm rev})$. 
\end{itemize}

Let us justify this construction. First of all, note that $\omega(A_n)=-w_\circ w(A_\circ)$. Secondly, it is easy to show that the alcove indexed by $\overline{q}=0$ in the gallery $\gamma(K)$ associated to some admissible subset $K$ in ${\mathcal A}(\Gamma^{\rm rev})$ is $-\ceil{\kappa_0(K)}(A_\circ)$; indeed, this is true for $K=\emptyset$, so, for an arbitrary $K$, one only needs to apply $\kappa_0(K)$ to the alcove indexed by $\overline{q}=0$ in $\gamma(\emptyset)$. We conclude that the alcove indexed by $\overline{q}=0$ in $J^{\rm rev}$ is $\omega(A_n)$ since $\ceil{\floor{w_\circ w}}=w_\circ w$. This means that $\gamma^\omega$ is a gallery. 

\begin{proposition}\label{admgal} The gallery $\gamma^\omega$ coincides with $\gamma(J^{\rm rev})$. In particular, we have $\mu(J^{\rm rev})=w_\circ(\mu(J))$. 
\end{proposition}

\begin{proof}
We will show that the admissible foldings corresponding to the two galleries coincide. In other words, we will prove that $\Gamma(\gamma^\omega)=\Gamma(J^{\rm rev})$, cf. the notation in Subsection \ref{chainroots}. Let 
\[\Gamma(J)=(\{(\gamma_i,\gamma_i')\}_{i\in I},\gamma_\infty)\,,\;\;\;\;\Gamma(J^{\rm rev})=(\{(\delta_i,\delta_i')\}_{i\in I},\delta_\infty)\,,\;\;\mbox{ and }\;\;\Gamma(\gamma^\omega)=(\{(\varepsilon_i,\varepsilon_i')\}_{i\in I},\varepsilon_\infty)\,.\]
By definition, the initial segments in $\Gamma(J^{\rm rev})$ and $\Gamma(\gamma^\omega)$ corresponding to $i\in \overline{I}$ coincide. We will now show that $\delta_i=\varepsilon_i$, for all $i\in[n]$; similarly, it can be shown that $\delta_i'=\varepsilon_i'$ and $\delta_\infty=\varepsilon_\infty$. Assume that $k_t<i\le k_{t+1}$ for some $t$ in $\{0,1,\ldots,s\}$ (if $t=0$ or $t=s$, one of the two inequalities is missing). Based on definitions and the fact that $r_{k_p}'=w_\circ^\lambda r_{j_{s+1-p}}w_\circ^\lambda$, we have
\begin{equation}\label{lrfolds}\begin{array}{ll}\delta_i&\!\!\!\!=\floor{w_\circ w}r_{k_1}'\ldots r_{k_t}'(\beta_i')=w_\circ w r_{j_s}r_{j_{s-1}}\ldots r_{j_{s+1-t}}w_\circ^\lambda(w_\circ^\lambda(\beta_{n+1-i}))\\[0.05in]&\!\!\!\!=w_\circ w r_{j_s}r_{j_{s-1}}\ldots r_{j_{s+1-t}}(\beta_{n+1-i})=w_\circ ur_{j_1}\ldots r_{j_{s-t}}(\beta_{n+1-i})\,.\end{array}\end{equation}
On the other hand, note that $\varepsilon_i$ is determined by $A_{i-1}'=\omega(A_{n+1-i})$ and $F_i'=\omega(F_{n+1-i})$. More precisely, we have $\varepsilon_i=-w_\circ(-\gamma_{n+1-i}')$. The proof is completed by observing that $j_{s-t}\le n+1-i<j_{s+1-t}$, which implies that 
\begin{equation}\label{lrfolds1}\varepsilon_i=w_\circ(\gamma_{n+1-i}')=w_\circ ur_{j_1}\ldots r_{j_{s-t}}(\beta_{n+1-i})\,.\end{equation}
Indeed, the expressions for $\delta_i$ and $\varepsilon_i$ in (\ref{lrfolds}) and (\ref{lrfolds1}) coincide.
\end{proof}

\subsection{The map $J\mapsto J^{\rm rev}$ and root operators} We will now present the main result related to the map $J\mapsto J^{\rm rev}$, which involves its commutation with the root operators. In order to do this, we need two lemmas. We will use once again the notation from Subsection \ref{comb-adm-seq}. In particular, given $J$ in ${\mathcal A}(\Gamma)$ as above and a simple root $\alpha_p$, we consider the set $I(J,p)$ and the sequence $\Sigma(J,p)$. We let $\Gamma(J)=(\{(\gamma_i,\gamma_i')\}_{i\in I},\gamma_\infty)$ and 
\begin{equation}\label{ijpsjp}\begin{array}{l}
I(J,p)\cap\overline{I}=\{\overline{i}_1<\ldots<\overline{i}_c\}\,,\;\;\;\;I(J,p)\setminus\overline{I}=\{{i}_1<\ldots<{i}_d\}\,,\\[0.05in]\Sigma(J,p)=(\overline{\sigma}_1,\ldots,\overline{\sigma}_{c},{\sigma}_1,\ldots,{\sigma}_{d},\sigma_{d+1})\,,\end{array}\end{equation}
where $c,d\ge 0$. Also recall that we set 
\[L(J)=\{l_i\}_{i\in I}\,,\;\;\;\;L(\emptyset)=\{l_i^\emptyset\}_{i\in I}\,,\;\;\;\;l_\infty^p:=\langle\mu(J),\alpha_p^\vee\rangle\,,\;\;\;\; u:=\kappa_0(J)\,,\;\;\;\;w:=\kappa_1(J)\,.\]

\begin{lemma}\label{lem1} If $c>0$, we have $\overline{\sigma}_1=\ldots=\overline{\sigma}_{c-1}=(1,-1)$, and either $\overline{\sigma}_c=(1,-1)$ or $\overline{\sigma}_c=(1,1)$; in the first case we have $w_\circ^\lambda u^{-1}(\alpha_p)<0$, while in the second one we have $w_\circ^\lambda u^{-1}(\alpha_p)>0$. If $c=0$, then $w_\circ^\lambda u^{-1}(\alpha_p)<0$. 
\end{lemma}

\begin{proof} We start by noting that, for $i\in\overline{I}$, the hyperplane $H_{-|\gamma_i|,l_i}$ is obtained from the hyperplane $H_{-\beta_i,l_i^\emptyset}=H_{\beta_i,0}$ by applying a nonaffine reflection;  therefore, $l_i=0$ for $i\in\overline{I}$. By Propositions \ref{gap} and \ref{combprop}, we can have $l_{\overline{i}_1}=\ldots=l_{\overline{i}_c}=0$ only if $\overline{\sigma}_i$ is as above, for $i=1,\ldots,c$. 

Let $\overline{J}:=J\cap\overline{I}$, and 
\begin{align*}&\Gamma(\overline{J})=(\{(\overline{\gamma}_i,\overline{\gamma}_i')\}_{i\in I},\overline{\gamma}_\infty)\,,\;\;\;
L(\overline{J})=\{\overline{l}_i\}_{i\in I}\,,\;\;\;\overline{l}_\infty^p:=\langle\mu(\overline{J}),\alpha_p^\vee\rangle\,,\\&I(\overline{J},p)=\{\overline{i}_1<\ldots<\overline{i}_c<h_1<\ldots<h_e\}\,,\;\;\;\; \Sigma(\overline{J},p)=(\overline{\sigma}_1,\ldots,\overline{\sigma}_{c},\pi_1,\ldots,\pi_e,\pi_{e+1})\,,
\end{align*}
where $e\ge 0$. Let $\beta:=|u^{-1}(\alpha_p)|$, and $\sigma:={\rm sgn}(u^{-1}(\alpha_p))={\rm sgn}(u(\beta))$. 

Assume first that $\langle\lambda,\beta^\vee\rangle=|\langle u(\lambda),\alpha_p^\vee\rangle|=|\langle\mu(\overline{J}),\alpha_p^\vee\rangle|=|\overline{l}_\infty^p|$ is nonzero. For $i=1,\ldots,e$, we have $\beta_{h_i}=\beta$ and $l_{h_i}^\emptyset=i$, which implies $\overline{l}_{h_i}=\sigma i$ and $\overline{l}_\infty^p=\sigma(e+1)$. If $\sigma=1$, we must have  $c>0$ and $\overline{\sigma}_c=(1,1)$ (by Proposition \ref{gap}). 
Similarly, if $\sigma=-1$, we must have $c>0$ (for this we also need Proposition \ref{combprop} (2)) and $\overline{\sigma}_c=(1,-1)$. 
Finally, the root $u^{-1}(\alpha_p)=\sigma\beta$ does not belong to the parabolic subroot system $\Phi_\lambda$ corresponding to $W_\lambda$, so $\sigma={\rm sgn}(u^{-1}(\alpha_p))={\rm sgn}(w_\circ^\lambda u^{-1}(\alpha_p))$. Indeed, if $\delta$ is a simple root in $\Phi_\lambda$, then $s_\delta$ sends a root in $\Phi^+\setminus\Phi_\lambda$ to another such root.

Now assume that $\langle\lambda,\beta^\vee\rangle=\overline{l}_\infty^p=0$, in which case we necessarily have $e=0$. If $c>0$, we must have  $\overline{\sigma}_c=(1,-1)$ and $\pi_1=1$ (by Propositions \ref{gap} and \ref{combprop}). But $\pi_1={\rm sgn}(\langle\overline{\gamma}_\infty,\alpha_p^\vee\rangle)={\rm sgn}(\langle u(\rho),\alpha_p^\vee\rangle)={\rm sgn}(u^{-1}(\alpha_p))$. On the other hand, $u^{-1}(\alpha_p)=\beta$ lies in $\Phi_\lambda$, so $w_\circ^\lambda u^{-1}(\alpha_p)<0$. If $c=0$, we must have $\pi_1=1$. The case $c=0$ is completely similar. 
\end{proof}

In addition to the notation in (\ref{ijpsjp}) related to the admissible subset $J$, we need the following one related to $J^{\rm rev}$:
\begin{equation}\label{ijrpsjrp}\begin{array}{l}
I(J^{\rm rev},p^*)\cap\overline{I}=\{\overline{h}_1<\ldots<\overline{h}_e\}\,,\;\;\;\;I(J^{\rm rev},p^*)\setminus\overline{I}=\{{h}_1<\ldots<{h}_f\}\,,\\[0.05in]\Sigma(J^{\rm rev},p^*)=(\overline{\pi}_1,\ldots,\overline{\pi}_{e},{\pi}_1,\ldots,{\pi}_{f},\pi_{f+1})\,.\end{array}\end{equation}
Let us define $\sigma_0\in\{-1,1\}$ by $\overline{\sigma}_c=(1,\sigma_0)$ if $c>0$, and by $\sigma_0:=-1$ if $c=0$. We define $\pi_0$ similarly, based on $e$ and $\overline{\pi}_e$. 
Given a pair of integers $(a,b)$, we also set $-(a,b):=(-a,-b)$. 

\begin{lemma}\label{lem2} We have $d=f$, as well as $h_j=n+1-i_{d+1-j}$ for $j=1,\ldots,d$ and $\pi_j=-\sigma_{d+1-j}$ for $j=0,1,\ldots,d+1$.
\end{lemma}

\begin{proof}
Let $\Gamma(J^{\rm rev})=(\{(\delta_i,\delta_i')\}_{i\in I},\delta_\infty)$. We will show that $\gamma_i=\pm\alpha_p$ implies $\delta_{n+1-i}'=\mp\alpha_{p^*}$, where $i\in[n]$; this, in turn, immediately implies $d=f$, as well as $h_j=n+1-i_{d+1-j}$ and $\pi_j=-\sigma_{d+1-j}$ for $j=1,\ldots,d$. Indeed, recall the setup related to the definition of $J^{\rm rev}$ in (\ref{defj})-(\ref{defjrev}), and assume that $k_t\le n+1-i< k_{t+1}$ for some $t$ in $\{0,1,\ldots,s\}$ (if $t=0$ or $t=s$, one of the two inequalities is missing); by (\ref{lrfolds}), we have
\[ \delta_{n+1-i}'=\floor{w_\circ w}r_{k_1}'\ldots r_{k_t}'(\beta_{n+1-i}')=w_\circ ur_{j_1}\ldots r_{j_{s-t}}(\beta_{i})=w_\circ(\gamma_i)\,,\]
where the last equality follows from the fact that $j_{s-t}<i\le j_{s+1-t}$. 

At this point, it suffices to show that $\pi_{d+1}=-\sigma_0$. By Proposition \ref{jrev} and an easy computation, we have 
\[\pi_{d+1}={\rm sgn}(\langle w(J^{\rm rev})(\rho),\alpha_{p^*}\rangle)={\rm sgn}(\langle w_\circ u w_\circ^\lambda(\rho),\alpha_{p^*}\rangle)=-{\rm sgn}(w_\circ^\lambda u^{-1}(\alpha_p))=-\sigma_0\,;\]
the last equality is the content of Lemma \ref{lem1}.
\end{proof}

\begin{theorem}\label{comm-f-rev} A root operator $F_p$ is defined on the admissible subset $J$ if and only if $E_{p^*}$ is defined on $J^{\rm rev}$, and we have
\[(F_p(J))^{\rm rev}=E_{p^*}(J^{\rm rev})\,.\]
\end{theorem}

\begin{proof} We use the setup above, particularly (\ref{defj})-(\ref{defjrev}) and (\ref{ijpsjp})-(\ref{ijrpsjrp}). We will compare $F_p(J)$ and $E_{p^*}(J^{\rm rev})$ in several cases. Let us assume first that $c,e>0$. Consider the functions $f\::\:[-c+\frac{1}{2},d+1]\rightarrow {\mathbb R}$ and $g\::\:[-e+\frac{1}{2},d+1]\rightarrow {\mathbb R}$ defined by
\[f(x)=g_{J,p}\!\left(x+c-\frac{1}{2}\right)\,,\;\;\;\;g(x)=g_{J^{\rm rev},p^*}\!\!\left(x+e-\frac{1}{2}\right)\,.\]
Based on Proposition \ref{gap}, these functions and the following observations related to them will be used below (sometimes implicitly) in order to construct $F_p(J)$ and $E_{p^*}(J^{\rm rev})$. By Lemmas \ref{lem1} and \ref{lem2}, we have $f(0)=g(0)=0$ and $g'(x)=-f'(d+1-x)$, for all $x\in[0,d+1]\setminus\frac{1}{2}{\mathbb Z}$. This means that 
$x_0$ is the first global maximum of $f$ on $[0,d+1]$ if and only if $d+1-x_0$ is the last global maximum of $g$ on $[0,d+1]$. By Proposition \ref{combprop}, the local maxima of $f$ and $g$ can only be attained at integer points.

{\bf Case 0:} $F_p$ is not defined on $J$, so $M(J,p)=0$. We have $f(x)\le 0=f(0)$ for all $x$ in its domain, and therefore $g(x)\le g(d+1)$. This means that $E_{p^*}$ is not defined on $J^{\rm rev}$. 

In fact, the above reasoning allows us to prove that $F_p$ is defined on $J$ if and only if $E_{p^*}$ is defined on $J^{\rm rev}$. The remaining cases deal with this situation.

{\bf Case 1:} $m_F(J,p)\ne\infty$ and $k_F(J,p)\in I\setminus\overline{I}$. This case is illustrated by the example in Figure~\ref{root1} below, where the graph on the left is of the function $f$, while the one on the right is of the function $g$; the dashed lines show the effect of applying the root operators $F_p$ to $J$ and $E_{p^*}$ to $J^{\rm rev}$. Let $m_F(J,p)=i_j\in J$ and $k_F(J,p)=i_{j-1}\not\in J$, where $1<j\le d$. Using Lemma \ref{lem2} and the above observations, we have 
\[k_E(J^{\rm rev},p^*)=n+1-i_j=h_{d+1-j}\in J^{\rm rev}\;\;\;\;\mbox{and}\;\;\;\;m_E(J^{\rm rev},p^*)=h_{d+2-j}=n+1-i_{j-1}\not\in J^{\rm rev}\,.\] 
Hence, we have
\[(F_p(J))^{\rm rev}=((J\setminus\{i_j\})\cup\{i_{j-1}\})^{\rm rev}\;\;\;\;\mbox{and}\;\;\;\;E_{p^*}(J^{\rm rev})=(J^{\rm rev}\setminus\{n+1-i_j\})\cup\{n+1-i_{j-1}\}\,.\]
In order to prove that these two admissible subsets above coincide, it suffices to show that their intersections with $\overline{I}$ coincide. The second intersection is $J^{\rm rev}\cap\overline{I}$, while $(F_p(J))^{\rm rev}\cap\overline{I}$ is computed based on $w(F_p(J))$. But this computation is the same as the one leading to $J^{\rm rev}\cap\overline{I}$, because we have $w(F_p(J))=w(J)$  by Proposition \ref{propf} (3). 
\begin{figure}[ht]
\mbox{\epsfig{file=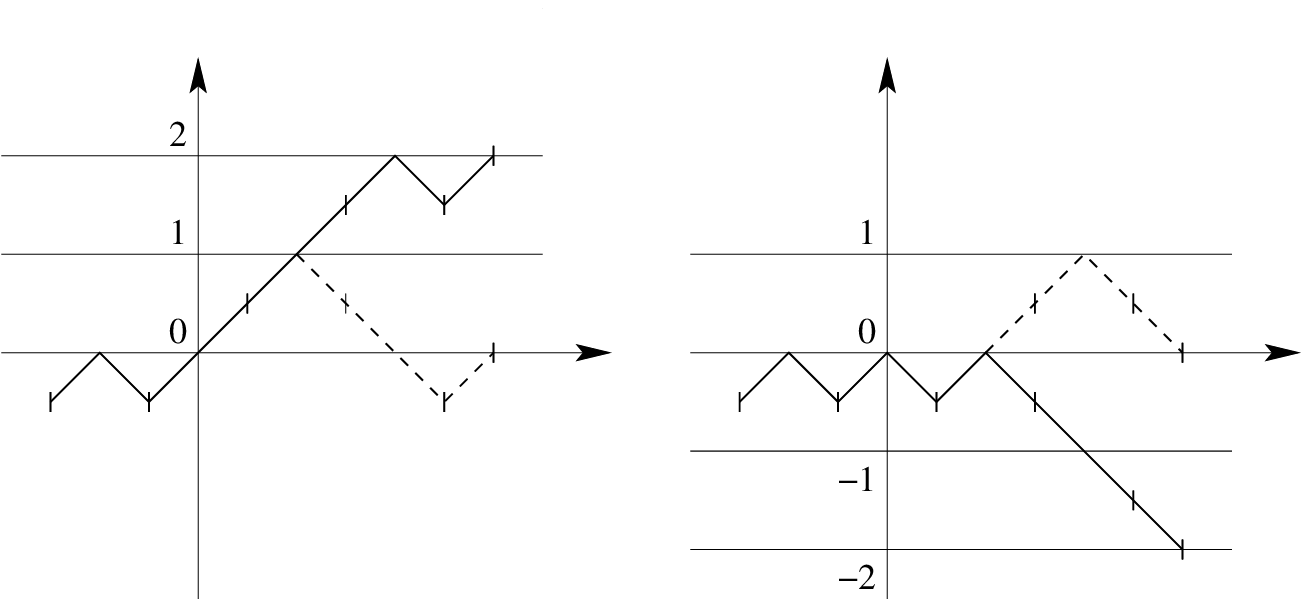}}
\caption{Case 1.}
\label{root1}
\end{figure}

{\bf Case 2:} $m_F(J,p)=\infty$ and $k_F(J,p)\in I\setminus\overline{I}$. This case is illustrated by the example in Figure~\ref{root2} below. In this case, the function $f$ has a unique global maximum (on its domain) at $d+1$, while $g$ has a unique global maximum on $[0,d+1]$ at 0, and $\overline{\pi}_e=(1,-1)$. Hence
\[k_F(J,p)=i_d\not\in J\,,\;\;\;\;k_E(J^{\rm rev},p^*)=\overline{h}_e\in J^{\rm rev}\,,\;\;\;\;m_E(J^{\rm rev},p^*)=h_1=n+1-i_d\not\in J^{\rm rev}\,.\]
Thus, we have
\[F_p(J)=J\cup\{i_d\}\;\;\;\;\mbox{and}\;\;\;\;(E_{p^*}(J^{\rm rev}))^{\rm rev}=((J^{\rm rev}\setminus\{\overline{h}_e\})\cup\{n+1-i_{d}\})^{\rm rev}\,.\]
These two admissible subsets coincide by a similar argument to the one used in Case 1. Note that we now need to use Proposition \ref{prope} (3) and Proposition \ref{jrev}, namely the fact that $w(E_{p^*}(J^{\rm rev}))=w(J^{\rm rev})=w_\circ \kappa_0(J) w_\circ^\lambda$. Indeed, this implies that the two admissible subsets above have the same initial key.
\begin{figure}[ht]
\mbox{\epsfig{file=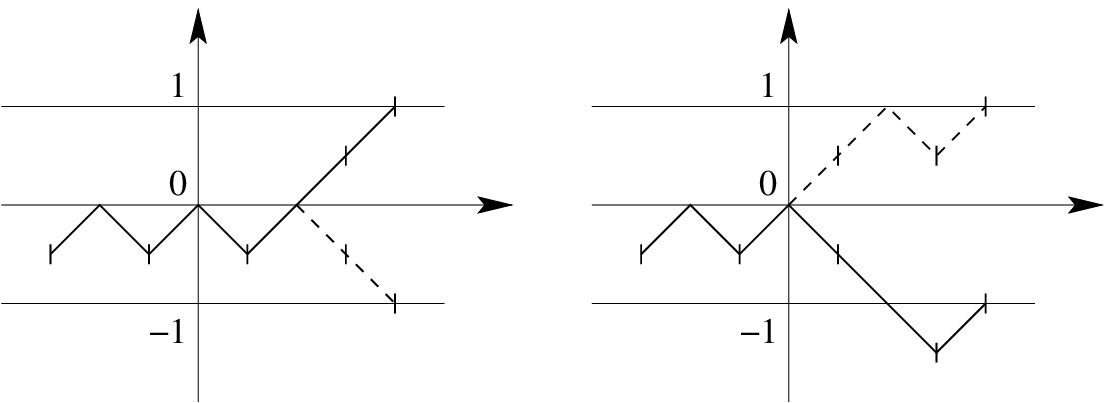}}
\caption{Case 2.}
\label{root2}
\end{figure}

{\bf Case 3:} $m_F(J,p)\ne\infty$ and $k_F(J,p)\in \overline{I}$. This case is illustrated by the example in Figure~\ref{root3} below. In this case, the function $f$ has its first global maximum (on its domain) at $1$, while $g$ has its last global maximum at $d$, and $\overline{\sigma}_c=(1,1)$. Hence
\[m_F(J,p)=i_1\in J\,,\;\;\;k_F(J,p)=\overline{i}_c\not\in J\,,\;\;\;k_E(J^{\rm rev},p^*)=n+1-i_1=h_d\in J^{\rm rev}\,,\;\;\;m_E(J^{\rm rev},p^*)=\infty.\]
Thus, we have
\[(F_p(J))^{\rm rev}=((J\setminus\{i_1\})\cup\{\overline{i}_c\})^{\rm rev}\;\;\;\;\mbox{and}\;\;\;\;E_{p^*}(J^{\rm rev})=J^{\rm rev}\setminus\{n+1-i_{1}\}\,.\]
These two admissible subsets coincide by a similar argument to the one used in Case 1.
\begin{figure}[ht]
\mbox{\epsfig{file=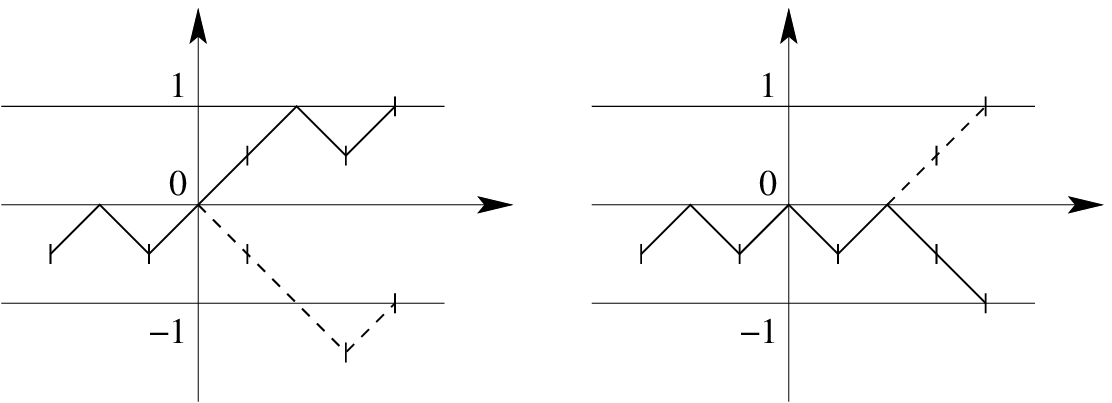}}
\caption{Case 3.}
\label{root3}
\end{figure}

{\bf Case 4:} $m_F(J,p)=\infty$ and $k_F(J,p)\in \overline{I}$. This case is illustrated by the example in Figure~\ref{root4} below. In this case we have $d=0$, $\overline{\sigma}_c=(1,1)$, $\sigma_1=1$, $\overline{\pi}_e=(1,-1)$, and $\pi_1=-1$. Hence
\[k_F(J,p)=\overline{i}_c\not\in J\,,\;\;\;\;k_E(J^{\rm rev},p^*)=\overline{h}_e\in J^{\rm rev}\,.\]
Thus, we have
\begin{equation}\label{fpeprev}(F_p(J))^{\rm rev}=(J\cup\{\overline{i}_c\})^{\rm rev}\;\;\;\;\mbox{and}\;\;\;\;E_{p^*}(J^{\rm rev})=J^{\rm rev}\setminus\{\overline{h}_e\}\,.\end{equation}
By Proposition \ref{propf} (3), we have $w(F_p(J))=s_pw(J)$, so
\[\kappa_0((F_p(J))^{\rm rev})=w_\circ s_pw(J)w_\circ^\lambda=s_{p^*}w_\circ w(J)w_\circ^\lambda=s_{p^*}\kappa_0(J^{\rm rev})\,.\]
By a completely similar proof to the one of Proposition \ref{prope} (3) in \cite{LP1}, we have $\kappa_0(J^{\rm rev}\setminus\{\overline{h}_e\})=s_{p^*}\kappa_0(J^{\rm rev})$. Therefore, we have $\kappa_0(E_{p^*}(J^{\rm rev}))=\kappa_0((F_p(J))^{\rm rev})$. This implies that the two admissible subsets in (\ref{fpeprev}) coincide.
\begin{figure}[ht]
\mbox{\epsfig{file=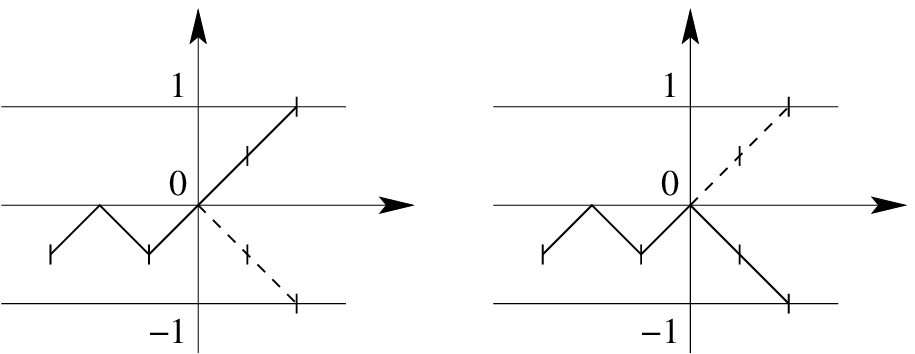}}
\caption{Case 4.}
\label{root4}
\end{figure}

We conclude the proof by discussing the case when $c=0$ or $e=0$. This is reduced to the simple observations below.
\begin{itemize}
\item If $c=0$ and $M(J,p)=0$ then $g_{J,p}$ attains its global maximum at $\frac{1}{2}$ (by Proposition \ref{combprop} (S2)), and therefore $g_{J^{\rm rev},p^*}$ attains its global maximum at $e+d-\frac{1}{2}$ and at the endpoint $e+d+\frac{1}{2}$; indeed, $\pi_{d+1}=1$ by Lemma \ref{lem2}. Hence, $E_{p^*}$ is not defined of $J^{\rm rev}$.
\item Case 1 is treated in the same way if $c=0$ or $e=0$. 
\item In Case 2 we cannot have $e=0$ because then $\pi_1=(-1,-1)$, and this is impossible by Proposition \ref{combprop} (S2). If $c=0$, then Case 2 is treated in the same way.
\item Case 3 does not make sense for $c=0$, and is treated in the same way if $e=0$. 
\item Case 4 does not exist.
\end{itemize} \end{proof}

We can summarize the construction in this section (based on Propositions \ref{gr-lamch} and \ref{jrev}) as follows: given the $\lambda$-chain $\Gamma$ (for a fixed dominant weight $\lambda$), we defined the $\lambda$-chain $\Gamma^{\rm rev}$, and given $J\in{\mathcal A}(\Gamma)$, we defined $J^{\rm rev}\in{\mathcal A}(\Gamma^{\rm rev})$. Hence we can map $J^{\rm rev}$ to an admissible subset $J^*\in{\mathcal A}(\Gamma)$ using Yang-Baxter moves, as it is described in Section \ref{yb-moves} and it is recalled below. To be more precise, let $R\::\:{\mathcal A}(\Gamma)\rightarrow {\mathcal A}(\Gamma^{\rm rev})$ denote the bijection $J\mapsto J^{\rm rev}$. On the other hand, we know from Proposition \ref{reversing} (2) that the $\lambda$-chains $\Gamma^{\rm rev}$ and $\Gamma$ can be related by a sequence of $\lambda$-chains $\Gamma^{\rm rev}=\Gamma_0,\,\Gamma_1,\,\ldots,\, \Gamma_m=\Gamma$ to which correspond Yang-Baxter moves $Y_1,\,\ldots,\,Y_m$. By Corollary \ref{isom-graphs}, the composition $Y:=Y_m\ldots Y_1$ does not depend on the sequence of intermediate $\lambda$-chains, and it defines a bijection from ${\mathcal A}(\Gamma^{\rm rev})$ to ${\mathcal A}(\Gamma)$. We let $J^*:=YR(J)$ and conclude that it is a bijection on ${\mathcal A}(\Gamma)$. The main result of this section, namely Theorem \ref{sch-inv-desc} below, now follows directly from Theorems \ref{comm-f-y} and \ref{comm-f-rev}.

\begin{remark}\label{no-move-init}
We claim that we can choose the $\lambda$-chains $\Gamma^{\rm rev}=\Gamma_0,\,\Gamma_1,\,\ldots,\, \Gamma_m=\Gamma$ such that their initial segments indexed by $\overline{1},\ldots,\overline{q}$ are identical. Indeed, this is true for $\Gamma^{\rm rev}$ and $\Gamma$ by definition. On the other hand, let us recall the correspondence between $\lambda$-chains and reduced words for the affine Weyl group element $v_{-\lambda}$ mentioned in Definition  \ref{def:lambda-chain}; most importantly, we recall from the proof of \cite[Lemma 9.3]{LP} that the moves $\Gamma_{i-1}\rightarrow\Gamma_{i}$ (for $i=1,\ldots,m$) correspond to Coxeter moves (on the mentioned reduced words) in this context. The claim is now justified by noting that two reduced words for $v_{-\lambda}$ with identical initial segments can be related by Coxeter moves which do not involve the mentioned initial segment. 
\end{remark}

\begin{theorem}\label{sch-inv-desc} The bijection $J\mapsto J^*$ constructed above coincides with Lusztig's involution $\eta_\lambda$ on the canonical basis. In other words, a root operator $F_p$ is defined on the admissible subset $J$ if and only if $E_{p^*}$ is defined on $J^*$, and we have 
\begin{equation}\label{cond1}(J_{\min})^*=J_{\max}\,,\;\;\;(J_{\max})^*=J_{\min}\,,\;\;\;\mbox{and}\;\;\;(F_p(J))^*=E_{p^*}(J^*)\,,\;\:\mbox{for}\;\:p=1,\ldots,r\,.\end{equation}
In particular, the map $J\mapsto J^*$ expresses combinatorially the self-duality of the poset ${\mathcal A}(\Gamma)$.
\end{theorem}

\begin{proof}
The first equality in (\ref{cond1}) is obvious. For the second one, note that $(J_{\max})^{\rm rev}=J_{\min}$ and use Remark \ref{no-move-init} to show that the map $Y$ fixes $J_{\min}$. The last equality follows directly from Theorems \ref{comm-f-y} and \ref{comm-f-rev}. Now recall that, based on the (directed colored graph) isomorphism in Corollary \ref{all-cryst}, we identified the vertex sets ${\mathcal B}_\lambda$ and ${\mathcal A}(\Gamma)$ of the corresponding directed colored graphs. By comparing (\ref{phipsieta0})-(\ref{phipsieta3}) with (\ref{cond1}), and by noting that the bijection specified by these conditions is unique, we conclude that the bijection $J\mapsto J^*$  coincides with $\eta_\lambda$ (via the isomorphism mentioned above). 
\end{proof}

\begin{remark} The above construction is analogous to the definition of Sch\"utzenberger's evacuation map (see, for instance, \cite{fulyt}). Below, we recall from Subsection \ref{cryst-inv} the three-step procedure defining this map and we discuss the analogy with our construction in the case of each step. 
\begin{enumerate}
\item REVERSE: We rotate a given semistandard Young tableau by 180$^\circ$. This corresponds to reversing its word, in the same way as we reversed the direction of our gallery, cf. Proposition \ref{admgal}. 
\item COMPLEMENT: We complement each entry via the map $i\mapsto w_\circ(i)$, where $w_\circ$ is the longest element in the corresponding symmetric group. This corresponds to using $w_\circ$ for the arbitrary Weyl group in the definition (\ref{defjrev}) of $J^{\rm rev}$.
\item SLIDE: We apply jeu de taquin on the obtained skew tableau. This corresponds to the Yang-Baxter moves $Y_1,\,\ldots,\,Y_m$ discussed above.
\end{enumerate}
\end{remark}

\begin{example}\label{ex-inv} Consider the Lie algebra $\mathfrak{sl}_3$ of type $A_2$, cf. Example \ref{simpleex}. Consider the dominant weight $\lambda=4\varepsilon_1+2\varepsilon_2$ and the following $\lambda$-chain:
\[
\begin{array}{ccccccccc}
 & \overline{1} & \overline{2}& \overline{3}& 1& 2& 3& 4& 5\\
\Gamma=(\!\!\!\!\!&\alpha_{12},&\alpha_{13},&\alpha_{23},&\alpha_{13},&{\underline{\alpha_{12}}},&\alpha_{13},&\underline{\alpha_{23}},&\alpha_{13})\,.
\end{array}
\]
Here we indicated the index corresponding to each root, using the notation in Subsection \ref{rev-lam-ch-j}; more precisely, we have $I=\{\overline{1}<\overline{2}<\overline{3}<1<2<3<4<5\}$ and $\overline{I}=\{\overline{1}<\overline{2}<\overline{3}\}$. By the defining relation (\ref{defgrev}), we have
\[
\begin{array}{ccccccccc}
& \overline{1} & \overline{2}& \overline{3}& 1& 2& 3& 4& 5\\
\!\!\!\!\Gamma^{\rm rev}\!=(\!\!\!\!\!&\underline{\alpha_{12}},&\alpha_{13},&\alpha_{23},&\alpha_{13},&\underline{\alpha_{23}},&\alpha_{13},&\underline{\alpha_{12}},&\alpha_{13})\,.
\end{array}
\]

Consider the admissible subset $J=\{2,\,4\}$. This is indicated above by the underlined roots in $\Gamma$. In order to define $J^{\rm rev}$, cf. (\ref{defjrev}), we need to compute
\[\kappa_0(J^{\rm rev})=w_\circ w(J)=(s_{12}s_{23}s_{12})(s_{12} s_{23})=s_{12}\,.\]
 Hence we have $J^{\rm rev}=\{\overline{1},\,2,\,4\}$. This is indicated above by the underlined positions in $\Gamma^{\rm rev}$. 

In order to transform the $\lambda$-chain $\Gamma^{\rm rev}$ into $\Gamma$, we need to perform a single Yang-Baxter move; this consists of reversing the order of the bracketed roots below:
\[
\begin{array}{lccccccccccc}
& &\overline{1} & \overline{2}& \overline{3}& 1& 2& 3& 4& 5\\
\Gamma^{\rm rev}&\!\!\!\!\!\!\!\!=(\!\!\!\!\!\!\!\!&\underline{\alpha_{12}},&\alpha_{13},&\alpha_{23},&\alpha_{13},&(\underline{\alpha_{23}},&\alpha_{13},&\underline{\alpha_{12}}),&\alpha_{13})&\longrightarrow
\\[0.1in]
 & &\overline{1} & \overline{2}& \overline{3}& 1& 2& 3& 4& 5\\
\Gamma&\!\!\!\!\!\!\!\!=(\!\!\!\!\!\!\!\!&\underline{\alpha_{12}},&\alpha_{13},&\alpha_{23},&\alpha_{13},&({{\alpha_{12}}},&\underline{\alpha_{13}},&\underline{\alpha_{23}}),&\alpha_{13})\,.
\end{array}
\]
The underlined roots indicate the way in which the Yang-Baxter move $J^{\rm rev}\mapsto Y(J^{\rm rev})=J^*$ works. All we need to know is that there are two saturated chains in Bruhat order between the permutations $u$ and $w$, cf. the notation in (\ref{setup1}):
\[u=s_{12}\lessdot s_{12}s_{23}\lessdot s_{12}s_{23}s_{12}=w\,,\;\;\;\;u=s_{12}\lessdot s_{12}s_{13}\lessdot s_{12}s_{13}s_{23}=w\,.\]
The first chain is retrieved as a subchain of $\Gamma^{\rm rev}$ and corresponds to $J^{\rm rev}$, while the second one is retrieved as a subchain of $\Gamma$ and corresponds to $J^{*}$. Hence we have $J^*=\{\overline{1},3,4\}$.
\end{example}

\section{Other Applications}\label{other-appl}

We can give an intrinsic explanation for the fact that the map $J\mapsto J^*$ is an involution on ${\mathcal A}(\Gamma)$; this explanation is only based on the results in Sections \ref{yb-moves} and \ref{sch-inv}, so it does not rely on Proposition \ref{phipsieta} (2). Let us first recall the bijections $R\::\:{\mathcal A}(\Gamma)\rightarrow {\mathcal A}(\Gamma^{\rm rev})$ and $Y\::\:{\mathcal A}(\Gamma^{\rm rev})\rightarrow {\mathcal A}(\Gamma)$ defined above. We claim that $YR=R^{-1}Y^{-1}$, which would prove that the composition $YR$ is an involution. In the same way as we proved Theorem \ref{sch-inv-desc} (that is, as a direct consequence of Theorems \ref{comm-f-y} and \ref{comm-f-rev}), we can verify that the composition $R^{-1}Y^{-1}$ satisfies the conditions in (\ref{cond1}). Since these conditions uniquely determine the corresponding map from ${\mathcal A}(\Gamma)$ to itself, our claim follows. 

\begin{remark} According to the above discussion, we have a second way of realizing Lusztig's involution $\eta_\lambda$ on the canonical basis, namely as $R^{-1}Y^{-1}$. In some sense, this is the analog of the construction of the evacuation map based on the {\em promotion} operation (see, for instance, \cite[p. 184]{fulyt}). To be more precise, the mentioned procedure has the following three steps.
\begin{enumerate}
\item Perform a sequence of sliding operations into the upper left corner of the given semistandard Young tableau, from which entries are removed successively. 
\item Place the removed entries into the corresponding outside corners that are vacated as a result of the sliding operations.
\item Complement the entries of the newly obtained filling of the corresponding Young diagram. 
\end{enumerate}
In one word, the sliding operations precede the complementation.
\end{remark}

We have the following corollary of Propositions \ref{jrev} and \ref{admgal}. According to this corollary, the alcove path model reveals an interesting feature of Lusztig's involution, which does not seem to be known even in type $A$. More precisely, it easily follows from our previous results that the involution $J\mapsto J^*$ interchanges the initial and the final keys in the sense mentioned below. 

\begin{corollary} \label{conjkeys}
For any $J\in {\mathcal A}(\Gamma)$, we have
\begin{equation}\label{cond2}\mu(J^*)=w_\circ(\mu(J))\,,\;\;\;\;\kappa_0(J^*)=\floor{w_\circ \kappa_1(J)}\,,\;\;\;\;\kappa_1(J^*)=\floor{w_\circ \kappa_0(J)}\,.\end{equation}
\end{corollary}

\begin{proof} The first equality follows directly from Proposition \ref{admgal} and the fact that a Yang-Baxter move preserves the weight of an admissible subset (cf. Theorem \ref{weightpres}). The second equality follows from the definition of $J^{\rm rev}$ in (\ref{defjrev}) combined with the fact that $\kappa_0(J^*)=\kappa_0(J^{\rm rev})$; the latter claim is a direct consequence of Remark \ref{no-move-init}. The third equality follows from (\ref{keys-rev}) and the fact that a Yang-Baxter move preserves the Weyl group element $w(\,\cdot\,)$ associated to an admissible subset (cf. (\ref{w-pres})).
\end{proof}


Recall the Demazure module $V_{\lambda,u}$ and its character $ch(V_{\lambda,u})$. Theorem \ref{charform} (2) provides a formula for this character. We now give a new formula, which we prove by setting up a bijection between the combinatorial objects indexing its terms and the combinatorial objects corresponding to the formula in Theorem \ref{charform} (2). 

\begin{theorem}\label{demform}  For any $u\in W$ and any $\lambda$-chain $\Gamma$, we have
\[ch(V_{\lambda,u})=\stacksum{J\in{\mathcal A}(\Gamma)}{w(J)\le u} e^{\mu(J)}\,.\]
\end{theorem}

\begin{proof}
By (\ref{w-pres}) and Theorem \ref{weightpres}, it suffices to consider a $\lambda$-chain $\Gamma$ (and the corresponding index set $I$) having the special form discussed at the beginning of Section \ref{sch-inv}. Let us assume first that $u$ is a maximal (left) coset representative modulo $W_\lambda$. We know from Theorem \ref{charform} (2) that
\begin{equation}\label{demchar}ch(V_{\lambda,u})=\sum e^{u r_{\overline{j}_1}\ldots r_{\overline{j}_a}\widehat{r}_{j_1}\ldots\widehat{r}_{j_s}(\lambda)}\,,\end{equation}
where the summation is over all subsets $J=\{\overline{j}_1<\ldots<\overline{j}_a<j_1<\ldots<j_s\}$ of $I$ such that we have a saturated decreasing chain in Bruhat order 
\[u\gtrdot ur_{\overline{j}_1}\gtrdot \ldots\gtrdot ur_{\overline{j}_1}\ldots r_{\overline{j}_a}\gtrdot ur_{\overline{j}_1}\ldots r_{\overline{j}_a}r_{j_1}\gtrdot\ldots\gtrdot ur_{\overline{j}_1}\ldots r_{\overline{j}_a}r_{j_1}\ldots r_{j_s}\,;\]
here it is assumed that $J\cap\overline{I}=\{\overline{j}_1<\ldots<\overline{j}_a\}$. Let $u':=w_\circ u r_{\overline{j}_1}\ldots r_{\overline{j}_a}$, which is a minimal coset representative modulo $W_\lambda$. There is a unique subset $\{\overline{k}_1<\ldots<\overline{k}_b\}$ of $\overline{I}$ such that 
\[1\lessdot r_{\overline{k}_1}\lessdot r_{\overline{k}_1}r_{\overline{k}_2}\lessdot\ldots\lessdot r_{\overline{k}_1}\ldots r_{\overline{k}_b}=u'\]
is a saturated increasing chain in Bruhat order from 1 to $u'$ (cf. Dyer \cite{Dyer}). Thus, $K:=\{\overline{k}_1<\ldots<\overline{k}_b<j_{1}<\ldots<j_s\}$ is an admissible subset. In fact, the map $J\mapsto K$ is a bijection between the subsets $J$ in (\ref{demchar}) and the admissible subsets $K$ with $\kappa_0(K)\ge w_\circ u$. Hence we have
\[ch(V_{\lambda,u})=\sum e^{w_\circ(\mu(K))}=\sum e^{\mu(K^*)}\,,\]
where the summations are over all admissible subsets $K$ with $\kappa_0(K)\ge w_\circ u$. But, by Corollary \ref{conjkeys} and the properties of the Bruhat order summarized in \cite[Lemma 2.1]{deoscb}, the latter condition is equivalent to $\kappa_1(K^*)=w(K^*)\le \floor{u}$.  The theorem  now follows by using the fact that $ch(V_{\lambda,u})=ch(V_{\lambda,\ceil{u}})$, for any $u$ in $W$, as well as the equivalence of $w(J)\le u$ and $w(J)\le \floor{u}$, where $J$ is an admissible subset (cf. \cite[Lemma 2.1]{deoscb}).
\end{proof}

\begin{remark}
Theorem \ref{demform} is the analog of the Demazure character formula due to Littelmann \cite{Li1}, \cite[Theorem 9.1]{litcrp}. Compared to the Demazure character formula in Theorem \ref{charform} (2), the one above has the advantage of realizing {\em all} Demazure characters $ch(V_{\lambda,u})$ (for a fixed $\lambda$) in terms of the {\em same} combinatorial objects, i.e., in terms of certain subsets of ${\mathcal A}(\Gamma)$. 
\end{remark}


\begin{thebibliography}{10}

\bibitem{bazcbq}
A.~Berenstein and A.~Zelevinsky.
\newblock Canonical bases for the quantum group of type {$A\sb r$} and
  piecewise-linear combinatorics.
\newblock {\em Duke Math. J.}, 82:473--502, 1996.

\bibitem{baztpm}
A.~Berenstein and A.~Zelevinsky.
\newblock Tensor product multiplicities, canonical bases and totally positive varieties.
\newblock {\em Invent. Math.}, 143:77--128, 2001.

\bibitem{bfpmbo}
F.~Brenti, S.~Fomin, and A.~Postnikov.
\newblock Mixed {B}ruhat operators and {Y}ang-{B}axter equations for {W}eyl
  groups.
\newblock {\em Internat. Math. Res. Notices}, 8:419--441, 1999.

\bibitem{Cher}
I.~Cherednik.
\newblock Quantum {K}nizhnik-{Z}amolodchikov equations and affine root systems.
\newblock {\em Comm. Math. Phys.}, 150:109--136, 1992.

\bibitem{Dem}
M.~Demazure.
\newblock D\mbox{\'e}singularization des vari\mbox{\'et\'e}s de
  \mbox{S}chubert.
\newblock {\em Annales E.N.S.}, 6:53--88, 1974.

\bibitem{deoscb}
V.~V. Deodhar.
\newblock A splitting criterion for the {B}ruhat orderings on {C}oxeter groups.
\newblock {\em Comm. Algebra}, 15:1889--1894, 1987.

\bibitem{Dyer}
M.~J. Dyer.
\newblock Hecke algebras and shellings of {B}ruhat intervals.
\newblock {\em Compositio Math.}, 89(1):91--115, 1993.

\bibitem{fulyt}
W.~Fulton.
\newblock {\em Young {T}ableaux}, volume~35 of {\em London Math. Soc. Student
  Texts}.
\newblock Cambridge Univ. Press, Cambridge and New York, 1997.

\bibitem{GaLi}
S.~Gaussent and P.~Littelmann.
\newblock {LS}-galleries, the path model and {MV}-cycles.
\newblock {\em Duke Math. J.},  127:35--88, 2005.

\bibitem{HK}
A. Henriques and J. Kamnitzer.
\newblock Crystals and coboundary categories.
\newblock {\em Duke Math. J.},  132:191--216, 2006.

\bibitem{Hum}
J.~E.~Humphreys.
\newblock {\em Reflection {G}roups and {C}oxeter {G}roups}, volume~29 of {\em
  Cambridge Studies in Advanced Mathematics}.
\newblock Cambridge University Press, Cambridge, 1990.

\bibitem{kascqa}
M.~Kashiwara.
\newblock Crystalizing the $q$-analogue of universal enveloping algebras.
\newblock {\em Commun. Math. Phys.}, 133:249--260, 1990.

\bibitem{kascbm}
M.~Kashiwara.
\newblock Crystal bases of modified quantized enveloping algebra.
\newblock {\em Duke Math. J.}, 73:383--413, 1994.

\bibitem{kascb}
M.~Kashiwara.
\newblock On crystal bases.
\newblock In {\em Representations of {G}roups (Banff, AB, 1994)}, volume~16 of
  {\em CMS Conf. Proc.}, pages 155--197. Amer. Math. Soc., Providence, RI,
  1995.

\bibitem{Kost}
B.~Kostant.
\newblock Powers of the {E}uler product and commutative subalgebras of a
  complex simple {L}ie algebra.
\newblock{\em Invent. Math.}, 158:181--226, 2004.

\bibitem{LS1}
V.~Lakshmibai and C.~S. Seshadri.
\newblock Standard monomial theory.
\newblock In {\em Proceedings of the Hyderabad Conference on Algebraic Groups
  (Hyderabad, 1989)}, pages 279--322, Madras, 1991. Manoj Prakashan.

\bibitem{lasdcg}
A.~Lascoux.
\newblock Double crystal graphs.
\newblock In {\em Studies in {M}emory of {I}ssai {S}chur ({C}hevaleret/{R}ehovot, 2000)},
  volume 210 of {\em Progr. Math.}, pages 95--114. Birkh\"auser Boston, Boston,
  MA, 2003.

\bibitem{lasksb}
A.~Lascoux and M.-P. Sch\mbox{\"{u}}tzenberger.
\newblock Keys and standard bases.
\newblock In D.~Stanton, editor, {\em Invariant Theory and Tableaux}, volume~19
  of {\em The IMA Vol. in Math. and Its Appl.}, pages 125--144,
  Berlin-Heidelberg-New York, 1990. Springer-Verlag.

\bibitem{lecsc}
C.~Lecouvey.
\newblock {Schensted-type correspondence, plactic monoid, and jeu de taquin for type {$C\sb n$}}.
\newblock {\em J. Algebra}, 247:295--331, 2002.

\bibitem{lecsbd}
C.~Lecouvey.
\newblock {Schensted-type correspondences and plactic monoids for types {$B\sb n$} and {$D\sb n$}}.
\newblock {\em J. Algebraic Combin.}, 18:99--133, 2003.

\bibitem{leeajt}
M.~A.~A. van Leeuwen.
\newblock An analogue of jeu de taquin for {L}ittelmann's crystal paths.
\newblock {\em S\'em. Loth. Comb.}, 41, Art. B41b, 23 pp., 1998. 

\bibitem{LP}
C.~Lenart and A.~Postnikov.
\newblock Affine {W}eyl groups in {$K$}-theory and representation theory.
\newblock {\tt arXiv:math.RT/0309207}.

\bibitem{LP1}
C.~Lenart and A.~Postnikov.
\newblock A combinatorial model for crystals of {K}ac-{M}oody algebras.
\newblock {\tt arXiv:math.RT/0502147}, to appear in {\em Trans. Amer. Math. Soc.}

\bibitem{Li1}
P.~Littelmann.
\newblock {A Littlewood-Richardson rule for symmetrizable Kac-Moody algebras}.
\newblock {\em Invent. Math.}, 116:329--346, 1994.

\bibitem{Li2}
P.~Littelmann.
\newblock Paths and root operators in representation theory.
\newblock {\em Ann. of Math. (2)}, 142:499--525, 1995.

\bibitem{litcrp}
P.~Littelmann.
\newblock Characters of representations and paths in {${\mathfrak h}\sp\ast\sb
  {\mathbb R}$}.
\newblock In {\em Representation {T}heory and {A}utomorphic {F}orms (Edinburgh,
  1996)}, volume~61 of {\em Proc. Sympos. Pure Math.}, pages 29--49. Amer.
  Math. Soc., Providence, RI, 1997.

\bibitem{Li3}
P.~Littelmann.
\newblock {Contracting modules and standard monomial theory for symmetrizable {K}ac-{M}oody algebras}.
\newblock {\em J. Amer. Math. Soc.}, 11:551--567, 1998.

\bibitem{lltpm}
M. Lothaire.
\newblock The plactic monoid (by A. Lascoux, B. Leclerc, and J-Y. Thibon).
\newblock In {\em Algebraic Combinatorics on Words}, pages 144--172. Cambridge University Press, Cambridge, 2002.  

\bibitem{luscba}
G.~Lusztig.
\newblock Canonical bases arising from quantized enveloping algebras. {II}.
\newblock {\em Progr. Theoret. Phys. Suppl.}, 102:175--201, 1991.

\bibitem{lusiqg}
G.~Lusztig.
\newblock {\em Introduction to Quantum Groups}, volume~110 of {\em
  Progress in Mathematics}.
\newblock Birkh\"auser, Boston, 1993.

\bibitem{smgrgb}
S.~Morier-Genoud.
\newblock Rel{\`e}vement g{\'e}om{\'e}trique de la base canonique et involution
  de {S}ch{\"u}tzenberger. ({F}rench) [{G}eometrical lifting of the canonical
  base and {S}ch{\"u}tzenberger involution].
\newblock {\em C. R. Math. Acad. Sci. Paris}, 337:371--374, 2003.

\bibitem{shesjt}
J.~Sheats.
\newblock {A symplectic jeu de taquin bijection between the tableaux of {K}ing and of {D}e {C}oncini}.
\newblock {\em Trans. Amer. Math. Soc.}, 351:3569--3607, 1999.

\bibitem{St}
J.~R. Stembridge.
\newblock Combinatorial models for {W}eyl characters.
\newblock {\em Adv.\ Math.}, 168:96--131, 2002.


\end{thebibliography}
\end{document}